\documentclass[final,onefignum,onetabnum]{siamart190516}

\usepackage{lipsum}
\usepackage{amsfonts}
\usepackage{graphicx}
\usepackage{epstopdf}

\usepackage{amsmath}
\usepackage{subcaption}
\usepackage{url}
\usepackage{hyperref}
\usepackage{algorithm}
\usepackage{algpseudocode}
\usepackage{nicefrac}
\usepackage{microtype}
\usepackage{bm}
\usepackage{enumitem}
\usepackage{color}
\usepackage{xcolor}

\ifpdf
  \DeclareGraphicsExtensions{.eps,.pdf,.png,.jpg}
\else
  \DeclareGraphicsExtensions{.eps}
\fi

\newsiamremark{remark}{Remark}
\newsiamremark{hypothesis}{Hypothesis}
\crefname{hypothesis}{Hypothesis}{Hypotheses}
\newsiamthm{claim}{Claim}
\def\red{\textcolor{red}}
\def\blue{\textcolor{blue}}
\setitemize{noitemsep,topsep=0pt,parsep=0pt,partopsep=0pt}
\newcommand{\E}{\mathbb{E}}

\headers{AEGD: Adaptive Gradient Descent with Energy}{H. Liu and X. Tian}

\title{AEGD: Adaptive Gradient Descent with Energy\thanks{Submitted to the editors 5/31/2021.
\funding{We would like to acknowledge support for this project
from the National Science Foundation (NSF grant DMS-1812666).}}}

\author{Hailiang Liu\thanks{Department of Mathematics,
Iowa State University, Ames, IA 50011, USA
  (\email{hliu@iastate.edu}, \url{https://faculty.sites.iastate.edu/hliu/}).}
\and Xuping Tian\thanks{Department of Mathematics,
Iowa State University, Ames, IA 50011, USA 
  (\email{xupingt@iastate.edu}).}}

\usepackage{amsopn}

\ifpdf
\hypersetup{
  pdftitle={AEGD: Adaptive Gradient Descent with Energy},
  pdfauthor={H. Liu and X. Tian}
}
\fi

\begin{document}

\maketitle

\begin{abstract}
We propose AEGD, a new algorithm for first-order gradient-based optimization of non-convex objective functions, based on a dynamically updated ‘energy’ variable. 
The method is shown to be unconditionally energy stable, irrespective of the step size. We prove energy-dependent convergence rates of AEGD for both non-convex and convex objectives, which for a suitably small step size recovers desired convergence rates for the batch gradient descent. We also provide an energy-dependent bound on the stationary convergence of AEGD in the stochastic non-convex setting.
The method is straightforward to implement and requires little tuning of hyper-parameters. Experimental results demonstrate that AEGD works well for a large variety of optimization problems: it is robust with respect to initial data, capable of making rapid initial progress. The stochastic AEGD shows comparable and often better generalization performance than SGD with momentum for deep neural networks. The code is available at \url{https://github.com/txping/AEGD}.
\end{abstract}



\section{Introduction}
From a  mathematical perspective, training neural networks is a high-dimensional non-convex optimization problem, and the dynamics of the training process is incredibly complicated. Despite this, stochastic gradient descent (SGD) \cite{RM51} and its variants have proven to be extremely effective for training neural networks in practice (see, for example, \cite{BCN18, GBC16}). Such stochastic training often reduces to solving the following unconstrained minimization problem
\begin{equation}\label{main}
\min_{\theta \in \mathbb{R}^n} f(\theta)
\end{equation}
where $f: \mathbb{R}^n \to \mathbb{R}$ 
is a finite-sum function, defined as  
$
f(\theta):=\frac{1}{m} \sum_{i=1}^m f_i(\theta)$, 
where $f_i(\theta):= L(F(x_i; \theta); y_i)$ is the loss of a given machine learning (ML) model $F(\cdot; \theta)$ on the training data $\{x_i, y_i\}$, parametrized by $\theta$.  For many practical applications, $f(\theta)$ is highly non-convex, and $F(\cdot; \theta)$ is chosen among deep neural networks (DNNs), known for their superior performance across various tasks. These deep models are heavily over-parametrized and require large amounts of training data. Thus, both $m$ and the dimension $n$ of $\theta$
can scale up to millions or even billions. These complications pose serious computational challenges.

Gradient descent (GD) or its variants 
are the method of choice for 
solving (\ref{main}). There are two variations of gradient decent, which differ in how much data we use to compute the gradient of the objective function. They are fully gradient decent, aka batch GD, and SGD. For large dataset, SGD is usually much faster by performing a parameter update for each training example.  
The GD method dates back to Cauchy (1847), and it updates $\theta$ according to:
\begin{equation}\label{gd}
\theta_{k+1}=\theta_k -\eta \nabla f(\theta_k),
\end{equation}
where $\eta> 0$ is the step size (called the learning rate in ML). GD has advantages of easy implementation and being fast for well-conditioned and strongly convex objectives, independent of the dimension of the underlying problem \cite{NW06}.  However, 
GD typically limits the step size to ensure numerical stability.
The step limitation issue can be largely resolved by implicit updates \cite{TA17},
\begin{equation}\label{be}
\theta_{k+1}=\theta_k -\eta \nabla f(\theta_{k+1}).
\end{equation}
This can be rewritten in the form of the celebrated proximal point algorithm (PPA) \cite{Ro76}: 
\begin{equation}\label{ppa}
\theta_{k+1}={\rm argmin}_{\theta} \left\{ f(\theta) +\frac{1}{2\eta}\|\theta-\theta_k\|^2\right\}.
\end{equation}
Such method has the advantage of being monotonically decreasing, which is guaranteed for any step size $\eta>0$. However, for non-quadratic function $f$, one would have to solve it by iteration, for example by the backward Euler method \cite{YPOO17}. A natural question is:

{\sl Can we adapt the GD (\ref{gd}) and improve stability with arbitrary step size as in (\ref{be}) but still easy to implement?}

We answer this question by proposing AEGD, a new method for efficient optimization that only requires first-order gradients. At its core, an auxiliary energy variable is used to adjust the step size at each update. The name AEGD is derived from the adaptive GD method with energy. This method is motivated by the energy dissipation structure hidden in the GD method with momentum (GDM) \cite{Po64}. We would like to emphasize that our AEGD features both a global and an element-wise formulation, while by empirical evidence, the latter appears to perform better in most cases.  

\subsection{Main Contributions}
To our best knowledge, AEGD  is the first algorithm that uses the energy update to stabilize GD in the optimization of non-convex objective functions. Meanwhile we show that AEGD still feathers the expected convergence rates for first-order gradient-based methods including GD. 
Some major benefits of AEGD are summarized below:
\begin{itemize}
\item It is shown to be unconditionally energy stable,  irrespective of the step size.
\item It allows a larger base step size than GD.
\item It converges fast even for objective functions that have a large condition number.
\item It is robust with respect to initial data.
\item It has visible advantages over the GD algorithm with momentum.
\end{itemize}
Most of the above advantages are kept for the stochastic AEGD, which will be discussed in Section 4. In particular, our numerical results also show that the stochastic AEGD has the potential for training neural networks with better generalization, at least on some benchmark machine learning tasks.



\subsection{Further Related Work} A numerical method bearing a direct relation to AEGD is the invariant energy quadratization (IEQ) approach introduced by \cite{Y16} and \cite{ZWY17} for gradient flows in the form of partial differential equations (PDEs). The IEQ approach uses modified quadratic energy to construct linear, unconditionally energy stable schemes for solving time-dependent PDEs \cite{LY19, LY20, SXY19}. The AEGD algorithm here is new in that we use an auxiliary energy variable to obtain unconditionally energy stable optimization algorithms.

In the field of stochastic optimization, there is an extensive volume of research for designing algorithms to speed up convergence.  Among others, two common approaches to accelerate SGD are momentum \cite{Po64} and adaptive step size \cite{DHS11, TH12, Ze12}.  The latter makes use of the historical gradient to adapt the step size, with renowned works like RMSProp \cite{TH12}, Adagrad \cite{DHS11}, and Adadelta \cite{Ze12} (see also the overview article \cite{Ru16}). The integration of momentum and adaptive step size led to Adam \cite{KB15}, which is one of the optimization algorithms that are widely used by the deep learning community. Many further advances have improved Adam \cite{ CK19, Do16, LH18, LJ+19, LH19, LXLS19, RKK19}.

The name momentum stems from an analogy to momentum in physics. A simple batch GD with constant momentum such as the heavy-ball (HB) method \cite{Po64} is known to enjoy the convergence rate of $O(1/k)$ for convex smooth optimization. Remarkably, with an adaptive momentum the Nesterov accelerated gradient (NAG) \cite{Ne83, SBC14} has the convergence rate up to the optimal $O(1/k^2)$. Recent advances showed that NAG has other advantages such as speeding up escaping saddle points \cite{JNJ17}, accelerating SGD or GD in non-convex problems \cite{Q99, Zhu18}.  One can also improve generalization with SGD in training DNNs with scheduled restart techniques \cite{OC15, RD17, WN+20}. Both first and second order optimization methods can be designed via inertial systems \cite{AC+20, SBC14, CB+21}.

Another major bottleneck for fast convergence of SGD lies in its variance \cite{Bo12, ST13, SQ96}. Different ideas have been proposed to develop variance reduction algorithm. Some representative works are SAGA \cite{DBL14}, SCSG \cite{LJCJ17}, SVRG \cite{JZ13}, Laplacian smoothing \cite{SW+19, WZGO20}, and iterative averaging methods \cite{PJ92, ZCM15}.

For a gradient based method the geometric property of $f$ often affects the convergence rates. For the non-convex $f$ we consider an old condition originally introduced by Polyak \cite{Po63}, who showed that it is a sufficient condition for gradient descent to achieve a linear convergence rate; see \cite{KNS20} for a recent convergence study under this condition. Such condition is often called the Polyak-\L ojasiewicz (PL) inequality in the literature since it is a special case of the inequality introduced by \L ojasiewicz in 1963 \cite{Lo63}. The works on convergence using the \L ojasiewicz inequality include \cite{Lo84, AMA05}. A generalization of the PL property for non-smooth optimization is the KL inequality \cite{Ku98, Boet09}. The KL inequality has been successfully used to study the convergence of the proximal algorithms \cite{AB13}, and the proximal alternating and projection methods \cite{ABRS10}.

\subsection{Organization} In Section 2, we motivate and present the AEGD method for non-convex optimization 
in both global and element-wise version. Theoretical results, including stability, convergence, and convergence rates, are given in Section 3, with technical proofs deferred to Appendix A. In Section 4 we present the stochastic AEGD and some theoretical results. Section 5 provides some experimental results to show the advantages of AEGD. We end this paper with concluding remarks in Section 6.

\subsection{Notation}  Throughout this paper, we denote $\{1, \cdots, m\}$ by $[m]$ for integer $m$. For vectors and matrices we use $\|\cdot\|$ to denote the $l_2$-norm. For a function $f: \mathbb{R}^n \to \mathbb{R}$, we use $\nabla f$ and $\nabla^2 f$ to denote its gradient and Hessian, and $\theta^*$ to denote a local minimizer of $f$. We also use the notation $\partial_i:= \partial_{\theta_i}$.
In the algorithm description we use $z=x/y$ to denote element-wise division if $x$ and $y$ are both vectors of size $n$; $x\odot y$ is element-wise product, and $x^2=x\odot x$.

\section{AEGD and Algorithm}\label{ae}
It is known one can accelerate gradients for certain objective functions by momentum, say by the HB scheme \cite{Po64} of form
\begin{equation*}
\theta_{k+1}=\theta_k -\eta \nabla f(\theta_k)+\mu(\theta_k-\theta^{k-1}),
\end{equation*}
where $0<\mu<1$ is a constant. This upon passing limit $\eta \to 0$ and $\mu \to 1$  with $\alpha=\frac{\eta}{(1-\mu)^2}$ kept fixed leads to 
\begin{align}\label{2nd}
\alpha \ddot \theta +\dot \theta=-\nabla f, \quad \alpha>0. 
\end{align}
The faster convergence and reduced oscillation may be expected from the perspective that the system is  driven by the energy dissipation law $ \dot E=-|\dot \theta|^2 \leq 0$ with the energy $E(t)=f(\theta(t))+\frac{\alpha}{2}|\dot \theta(t)|^2$.  Motivated by this energy dissipation structure, we define $g(\theta) =\sqrt{f(\theta)+c}$ with $c$ chosen so that $f(\theta)+c>0$. Set $r = g(\theta)$, then $r^2$ plays the role of energy and
$$
\nabla f  = 2r \nabla g, \quad r=g(\theta).
$$
This transformation leads to an augmented system 
\begin{subequations}\label{gr+}
\begin{align}
&\dot \theta=-2r \nabla g(\theta),  \\
&\dot r = \nabla g \cdot \dot \theta.
\end{align}
\end{subequations}
It allows the dissipation of the quadratic energy
$
\frac{d}{dt}(r^2)=-|\dot \theta|^2
$
in place of the  dissipation of nonlinear loss function:  $ \frac{d}{dt} f(\theta) =-|\dot \theta|^2$. Note also that the above transformation with $c$ taken as constant is linked back to the first order gradient flow $\dot \theta =-\nabla f$, instead of the inertial equation (\ref{2nd}). Based on (\ref{gr+}) we introduce the following update rule:
\begin{subequations}\label{e1}
\begin{align}
& \theta_{k+1}= \theta_k- 2\eta r_{k+1}\nabla g(\theta_k),\quad k=0,1,2,\cdots, \\
& r_{k+1}-r_k  = \nabla g(\theta_k) \cdot (\theta_{k+1}-\theta_k).
\end{align}
\end{subequations}
This is actually a linear algorithm 
since 
\begin{subequations}\label{e1+}
\begin{align}
& r_{k+1} =\frac{r_k}{1+2\eta |\nabla g(\theta_k)|^2}, \quad r_0=\sqrt{f(\theta_0)+c},\\
& \theta_{k+1}=\theta_k-2\eta r_{k+1}\nabla g(\theta_k), \quad k=0,1,2,\cdots,
\end{align}
\end{subequations}
which are easy to implement. In order to allow the use of different learning rate in each coordinate, we also propose an element-wise AEGD, set forth as follows:
\begin{subequations}\label{ee1}
\begin{align}
& r_{k+1,i} =\frac{r_{k,i}}{1+2\eta(\partial_i g(\theta_k))^2}, \quad i\in[n],\quad r_{0,i}=\sqrt{f(\theta_0)+c},\\
& \theta_{k+1,i}=\theta_{k,i}-2\eta r_{k+1,i}\partial_i g(\theta_k). \quad k=0,1,2,\cdots
\end{align}
\end{subequations}
This also implies
\begin{align}\label{ri}
r_{k+1,i}-r_{k,i}  = \partial_i g(\theta_k)(\theta_{k+1,i}-\theta_{k,i}), \quad i\in [n].
\end{align}
\begin{remark}
The essential idea behind \eqref{gr+} is the so called invariant energy quadratizaton (IEQ) strategy, originally introduced for developing linear and unconditionally energy stable schemes for gradient flows in the form of partial  differential equations \cite{Y16, ZWY17}. One may wonder whether  $r=(f+c)^\alpha$ is  admissible for building similar algorithms for all $\alpha \in (0, 1)$. To see this, we write the corresponding augmented system as  
$$
\dot \theta=-\alpha^{-1} r^{1/\alpha -1} \nabla g(\theta), \quad
\dot r = \nabla g \cdot \dot \theta.
$$
A similar implicit- explicit discretization can yield a linear update for $r$ if and only if 
$
\frac{1}{\alpha}-1=1 \Leftrightarrow \alpha=\frac{1}{2}.  
$
This explains why energy quadratization  has been used in the PDE community to name such strategy.
\end{remark}

\section{Theoretical Results}\label{con}
In this section, we will show AEGD is unconditionally energy stable, irrespective of the step size $\eta$, and obtain different convergence rates for non-convex, convex, and strongly convex functions. We say $f$ is $\alpha$-strongly convex if for all $u$ and $v$ we have $f(u)\geq f(v)+\langle \nabla f(v), u-v \rangle +\alpha\|u-v\|^2/2$; $f$ is $L$-smooth if $\|\nabla^2 f(u)\|\leq L$ for any $u\in \mathbb{R}^n$. For convex  ($\alpha=0$) and $L$-smooth functions, it is known that the convergence rate $O(1/k)$ for GD of form (\ref{gd}) is guaranteed if $\eta \leq \frac{1}{L}$ \cite{Ne98}.
This may suggest very small step-sizes in practice, which is a condition that may be violated in more complicated scenarios. In contrast, AEGD is energy stable for any $\eta>0$ (see Theorem \ref{thm1}), and converges for a large range of  $\eta$ (see Theorem \ref{thm2}). 
\subsection{Unconditional energy stability} 
\begin{theorem}
\label{thm1} (Energy stability and convergence)  Consider
$ \min \{f(\theta), \; \theta \in \mathbb{R}^n\},$
where $f(\theta)$ is differentiable and bounded from below so that $f(\theta)+c>0$  for some $c>0$.
Then \\
(i) AEGD (\ref{e1+})  is unconditionally energy stable in the sense that for any step size $\eta>0$,
 \begin{equation}\label{re}
 r_{k+1}^2=r^2_k -(r_{k+1} -r_k)^2- \eta^{-1} \|\theta_{k+1}-\theta_k\|^2,
 \end{equation}
 $r_k$ is strictly decreasing and convergent with $r_k \to r^*$ as $k\to \infty$, and also
  \begin{equation}\label{rev}
 \lim_{k\to \infty} \|\theta_{k+1}-\theta_k\|=0, \quad \sum_{j=0}^\infty\|\theta_{j+1}-\theta_j\|^2 \leq \eta (r^2_0-(r^*)^2).
  \end{equation}
  (ii) AEGD (\ref{ee1})  is unconditionally energy stable in the sense that for any step size $\eta>0$,
   \begin{equation}\label{rei}
   r^2_{k+1,i}=r^2_{k,i} -(r_{k+1,i} -r_{k,i})^2- \eta^{-1} (\theta_{k+1,i}-\theta_{k,i})^2,\quad i\in [n],
   \end{equation}
   $r_{k,i}$ is strictly decreasing and convergent with $r_{k,i} \to r_i^*$ as $k\to \infty$, and also
     \begin{equation}\label{rev1}
   \lim_{k\to \infty} \|\theta_{k+1}-\theta_k\|=0, \quad \sum_{j=0}^\infty\|\theta_{j+1}-\theta_j\|^2 \leq \eta \sum_{i=1}^n (r^2_{0,i}
   -(r^*_i)^2).
    \end{equation}
\end{theorem}
\begin{proof}
(i)  We use the two equations in (\ref{e1}) to derive
$$
2r_{k+1}(r_{k+1}-r_k) = 2r_{k+1} \nabla g(\theta_k)\cdot(\theta_{k+1}- \theta_k)=-4\eta r_{k+1}\|\nabla g(\theta_k)\|^2=-
\frac{1}{\eta}\|\theta_{k+1}-\theta_k\|^2.
$$
Upon rewriting with $2b(b-a) = b^2-a^2 + (b-a)^2$  we obtain equality (\ref{re}). From this we see  that $ r^2_k$ is monotonically decreasing (also bounded below), therefore convergent; so does $r_k$ since $r_k \geq 0$. Summation of (\ref{re}) over $k$ from $0, 1, \cdots$ yields (\ref{rev}). The proof of (ii) is entirely similar.
\end{proof}
\begin{remark}
It is worth pointing out that this theorem does not require the function $f$ satisfy L-smoothness or convexity assumption. The result asserts that the unconditional energy stability featured by AEGD also implies convergence of $\{r_k\}_{k\geq 0}$ for any $\eta>0$, and  the sequence $\{\|\theta_{k+1}-\theta_k\|\}_{k \geq 0}$ converges to zero at a rate of at least $1/\sqrt{k}$. But this is not sufficient---at least in general---to prove convergence of the sequence $\{\theta_k\}_{k\geq 0}$, when no further information is available about this sequence.  
\end{remark}
\subsection{Convergence and convergence rates} 
To understand the convergence behavior of AEGD (\ref{e1+}), we reformulate it as 
\begin{equation}\label{eq:energy}
 \theta_{k+1}= \theta_k- \eta_k \nabla f(\theta_k), \quad \eta_k:=\eta \frac{r_{k+1}}{g(\theta_k)}.
\end{equation}
Note that using $L$-smoothness  of $f$ and (\ref{eq:energy}) we have
\begin{align*}
f(\theta_{k+1}) & \leq
f(\theta_k) -\left(\frac{1}{\eta_k} -\frac{L}{2}\right)\|\theta_{k+1}-\theta_k\|^2 < f(\theta_k),
\end{align*}
when $\eta_k$ gets smaller so that $\eta_k<2/L$. Since $r_k$ is decreasing, after finite number of iterations $\eta_j$ can be ensured (by choosing $\eta$ to be suitably small if necessary) to fall below $2/L$ under which $f(\theta_j)$ turns into a strictly decreasing sequence, hence convergent. To ensure the convergence of $\{\theta_k\}$, one needs further info on the  geometrical property of the objective function, such as convexity, or the Kurdyka-Lojasiewicz (KL) property. The KL property at a point describes how the objective function can be made sharp through a concave mapping near that point. This property characterizes a rich function class, and often considered as a structural assumption in general nonconvex optimization (see \cite{AB13} for a detailed account on this property and its applications). In the present work we restrict to a special case of KL, called Polyak-Lojasiewics (PL) property: $f$ is PL if there exists $\mu>0$ such that 
\begin{align}\label{PL}
\frac{1}{2}\|\nabla f(\theta)\|^2 \geq \mu(f(\theta)-f^*),
\end{align}
where we assume that $\{\theta, \min f(\theta)=f^*\}$ is not empty. One such example $f(x)=x^2+3sin^2(x)$ is non-convex,  yet satisfying the PL inequality with $f^*=0$ and $\mu=1/32$.

{In Theorem \ref{thm2} below, we present convergence rates of AEGD in three different cases: nonconvex with the PL property, convex and strongly convex.} 
{For general $\eta$, convergence rates can depend on the behavior of $r_k$, which is part of the solution to the AEGD algorithm. Therefore, we present a hybrid result of convergence rates combining both $k$ and $r_k$ (a posterior estimates). }

\begin{theorem} \label{thm2} (Convergence rates)
Suppose $f$ is differentiable and bounded from below.
Let $\theta_k$ be the $k$-th iterate generated by AEGD \eqref{eq:energy}, then   
\begin{equation*}
\frac{1}{k}\sum_{j=0}^{k-1}\|\nabla g(\theta_j)\|^2
\leq \frac{\sqrt{f(\theta_0)+c}}{2\eta k r_k}.
\end{equation*}
We have convergence rates in three distinct cases:
\\
(i) 
$f$ is PL and $L$-smooth with a minimizer $\theta^*$. If $\max_{k_0\leq j<k} \eta_j \leq 1/L$ for some $k_0\geq 0$, then $\{\theta_k\}$ is convergent. Moreover,
\begin{align}\label{ctheta}
    & \sum_{k=k_0}^\infty\|\theta_{k+1}-\theta_k\|\leq \frac{4}{\sqrt{2\mu}} \sqrt{f(\theta_{k_0})-f(\theta^*)},\\\notag
    & f(\theta_k)-f(\theta^*)\leq e^{-c_0(k-k_0)r_k}(f(\theta_{k_0})-f(\theta^*)), \quad c_0:=\frac{\mu \eta}{\sqrt{f(\theta_{k_0})+c}}.
\end{align}
(ii) $f$ is convex and $L$-smooth with a global minimizer $\theta^*$. If $\max_{k_0\leq j<k} \eta_j \leq 1/L$ for some $k_0\geq 0$,  then
\begin{equation*}
f(\theta_k)-f(\theta^*) \leq \frac{c_1\|\theta_{k_0} -\theta^*\|^2}{2(k-k_0)r_k}, \; c_1:=\frac{\sqrt{f(\theta_{k_0})+c}}{\eta}.
\end{equation*}
(iii) $f$ is $\alpha$-strongly convex and  $L$-smooth with the global minimizer $\theta^*$. Assume that $\max_{k_0\leq j<k} \eta_j \leq \frac{2}{\alpha+L}$ for some $k_0 \geq 0$,  then
\begin{equation*}
\|\theta_k-\theta^*\| \leq e^{-c_2 (k-k_0) r_k}\|\theta_{k_0}-\theta^*\|, \quad c_2:=\frac{\alpha \eta}{\sqrt{f(\theta_{k_0})+c}}.
\end{equation*}
\end{theorem}
See Appendix A.1 for the proof.

\begin{remark} 
The gradient estimate in Theorem \ref{thm2}  when using  $g(\theta)=\sqrt{f(\theta)+c}$ gives 
$$
\frac{1}{k}\sum_{j=0}^{k-1}\|\nabla f(\theta_j)\|^2 \leq
\frac{2(f(\theta_0)+c)^{\frac{1}{2}}(F_k+c)}{\eta kr_k},\quad F_k=max_{j<k}f(\theta_j).
$$
For this estimate, neither $L$-smoothness nor step size restriction is needed.  This is in sharp contrast to the classical result for GD (see pages 29-31 in \cite{Ne98}).
\end{remark}

\begin{remark}
Regarding the convergence rates, a series of remarks are in order. 
\begin{enumerate} 
\item The convergence rate in (i) may be extended to the case when $f$ features the more general KL property; see \cite{AB13, ABRS10} for relevant techniques in establishing  convergence rates for proximal and other gradient methods.
\item If $r^*>0$, which is the case to be shown in Lemma \ref{tau} for $\eta<\tau$, $r_k$ on the right of the three estimates in Theorem \ref{thm2} can all be replaced by $r^*(<r_k)$, hence (i) guarantees the linear convergence rate when $f$ is PL;  
and (ii) recovers the usual sublinear rate $O(\frac{1}{k})$;  while with strict convexity, (iii) gives the usual linear convergence rate.
\item 
It is clear that all results are valid as long as $r_k \to 0$ slower than $1/k$. 
Our numerical tests indicate that for each objective function $f$ there exists a sharp threshold $\tilde \eta$ in the sense that convergence is ensured if $\eta \leq \tilde\eta$. However, identifying $\tilde \eta$ appears to be a challenging task in theory. 
\end{enumerate}

\end{remark}

\begin{remark}
Regarding the base step size $\eta$, we make further remarks. 
\begin{enumerate}
\item Empirical evidence (see Figure \ref{fig:rosen} (c)) shows that there exists a threshold index $J$ so that $\eta_j $ turns to decrease for $j >J$. Hence the assumptions in  (i)-(iii) are reasonable and readily met. 
\item Our results suggest that decaying $\eta$ at a later stage (say $k\geq k_0$ for some $k_0$) to ensure the sufficient conditions in Theorem \ref{thm2} can help to achieve all-time good performance of the AEGD algorithm. 
This comment applies to the element-wise AEGD as well (see Theorem  \ref{thm3}).
\end{enumerate}
\end{remark}
\subsection{Behavior of the energy} 
Since $r_k$ is strictly decreasing and $g(\theta_k)$ is positive and bounded from below and above, their relative ratio $\eta_k$ essentially depends on the behavior of $r_k$.  On the other hand, as $\eta \to 0$, the numerical solution $(\theta_k, r_k)$ may be shown to converge to the solution of the ODE system 
\begin{align}\label{rr}
& \dot \theta=-2r\nabla g,  \quad \dot r=\nabla g\cdot \dot \theta
\end{align}
at $t_k=k\eta$, subject to initial data $\theta(0)=\theta_0,  r(0)=g(\theta_0)$.  
For this system the level set $r(t)-g(\theta(t))=0$ is invariant for all time. Hence for a fixed but suitably small $\eta$, starting from $(\theta_0, g(\theta_0))$, the limit of $(\theta_k, r_k)$ as $k \to \infty$ must be approaching to $(\theta^*, g(\theta^*))$, that is   
$$
r^*=g(\theta^*)>0. 
$$
A natural question is whether a threshold for the base step size $\eta$ can be identified so that we will still have $r^*>0$ for any $\eta <\tau$. Indeed, we are able to obtain a sufficient condition to ensure this.
\begin{lemma} Suppose $f$ is $L$-smooth,  bounded from below by $f(\theta^*)$ and we have $\max \|\nabla f(\theta)\|\leq G_\infty$, then  $g$ is $L_g$-smooth with 
$$
g \geq g(\theta^*)=\sqrt{f(\theta^*)+c}, \quad L_g=\frac{1}{2g(\theta^*)} \left( L+ \frac{G_\infty^2}{2g^2(\theta^*)}\right). 
$$
Consider AEGD (\ref{e1+}), then  
$$
r_k> r^* > g(\theta^*)(1-\eta/\tau), \quad \tau:=\frac{2g(\theta^*)}{L_g (g(\theta_0))^2}.
$$  
This implies that $r^*>0$ if $\eta \leq \tau$.  
\end{lemma}\label{tau}
\begin{proof} 

For any $x, y\in \{\theta_k\}_{j=0}^T$ we have 
\begin{align*}
\|\nabla g(x)-\nabla g(y)\|
&
= \frac{1}{2}\bigg\|\frac{\nabla f(x)(g(y)-g(x))}{g(x)g(y)} + \frac{\nabla f(x)-\nabla f(y)}{g(y)}\bigg\|\\
&\leq \frac{G_\infty}{2(g(\theta^*))^2}|g(y)-g(x)| + \frac{1}{2g(\theta^*)}\|\nabla f(x)-\nabla f(y)\|\leq L_g\|x-y\|.
\end{align*}
Such $L_g$-smoothness of $g$ implies that 
\begin{align*}
g(\theta_{j+1}) & \leq g(\theta_j)+\nabla g(\theta_j)\cdot (\theta_{j+1}-\theta_j) +\frac{L_g}{2}\|\theta_{j+1}-\theta_j\|^2 \\
& = g(\theta_j)+r_{j+1}-r_j +\frac{L_g}{2}\|\theta_{j+1}-\theta_j\|^2.
\end{align*}
Take a summation over $j$ from $0$ to $k-1$ so that 
\begin{align*}
g(\theta_k)-g(\theta_{0}) &\leq r_k-r_0 + \frac{L_g}{2}\sum_{j=0}^{k-1} \|\theta_{j+1}-\theta_j\|^2 \\
&=  r_k-r_0 + \frac{L_g \eta }{2} \sum_{j=0}^{k-1} \left( 
r^2_j-r^2_{j+1} -(r_{j+1}-r_j)^2\right) \\
& =  r_k-r_0 + \frac{L_g \eta }{2} \left( 
r^2_0-r^2_k -\sum_{j=0}^{k-1}(r_{j+1}-r_j)^2
\right),
\end{align*}
where (\ref{re}) was used. Using $r_0=g(\theta_0)$ and $g(\theta_k)\geq g(\theta^*)$, we have for any $k$, 
$$
 \frac{L_g \eta }{2} \sum_{j=0}^{k-1}(r_{j+1}-r_j)^2
 + \frac{L_g \eta }{2} r^2_k -r_k +g(\theta^*) - \frac{L_g \eta }{2} (g(\theta_0))^2 \leq 0. 
$$
Passing to the limit as $k \to \infty$, we also have
$$
 \frac{L_g \eta }{2} \sum_{j=0}^{\infty}(r_{j+1}-r_j)^2
+ \frac{L_g \eta }{2} (r^*)^2 -r^* +g(\theta^*) (1-\eta/\tau) \leq 0. 
$$
Hence 
$
r^* > g(\theta^*)(1-\eta/\tau).
$ 
\end{proof}
\begin{remark} 1. Note that $\eta < \tau$ is only a sufficient condition, not necessary for $r^*>0$ to surely happen. \\
2. One may take a suitably large $\eta$ to gain initial rapid progress, and adjust $\eta$ at a later stage at $k=k_0$. A similar argument shows that $r^*>0$ is ensured by $\eta<\tau_1$ for $k \geq k_0$,
\begin{equation}\label{ee}
\eta <\tau_1:=\frac{2(g(\theta^*) +r_{k_0}- g(\theta_{k_0}))}{L_g r^2_{k_0}}. 
\end{equation}
\end{remark} 
{As we expected, for small $\eta$, AEGD  features same convergence rates as GD does, since the rates are essentially depending on the local geometry of $f$ near $\theta^*$. }
\subsection{Convergence results for the element-wise AEGD} 
Similar results also hold for the element-wise AEGD (\ref{ee1}), although analysis is more involved.  With the notation
$$
\eta_{ij}:=\eta r_{j+1,i}/g(\theta_j), \; i\in [n], \quad j=0,1,2,\cdots,
$$
we now present the main result for (\ref{ee1}) in the following.
\begin{theorem}\label{thm3} (Convergence rates)
Suppose $f$ is differentiable and bounded from below.
Let $\theta_k$ be the $k$-th iterate generated by the AEGD (\ref{ee1}), then 
\begin{equation*}
\min_{j<k}(\partial_i g(\theta_j))^2 
\leq \frac{1}{k}\sum_{j=0}^{k-1}(\partial_i g(\theta_j))^2
\leq \cfrac{\sqrt{f(\theta_0)+c}}{2\eta k r_{k,i}},\quad i\in[n].
\end{equation*}
We have convergence rates in three distinct cases:
\\
(i) $f$ is PL and $L$-smooth with a minimizer $\theta^*$. If $\max_{k_0\leq j\leq k}\eta_{ij}\leq\frac{1}{L}$ for $i\in[n]$ and some $k_0\geq0$, then
\begin{align*}
& f(\theta_k)-f(\theta^*)\leq e^{-c_0(k-k_0)\min_ir_{k,i}}(f(\theta_{k_0})-f(\theta^*)), \quad c_0:=\frac{\mu \eta}{\sqrt{f(\theta_{k_0})+c}}.
\end{align*}
(ii) $f$ is convex and $L$-smooth with a global minimizer $\theta^*$. If $\max_{k_0\leq j<k}\eta_{ij}  \leq\frac{1}{L}$ for $i\in[n]$ and some $k_0\geq 0$, then
$$
f(\theta_k)-f(\theta^*) \leq \frac{c_1\max_{k_0\leq j<k}\|\theta_j -\theta^*\|^2}{(k-k_0) \min_i r_{k,i}}, \; c_1:=\frac{\sqrt{f(\theta_{k_0})+c}}{\eta}.
$$
(iii) $f$ is $\alpha-$strongly convex and  $L$-smooth with the global minimizer $\theta^*$. Assume that
 $ \max_{k_0\leq j<k}\eta_{ij}  \leq \frac{2}{\alpha+L}$ for $i\in[n]$ and some $k_0\geq 0$,
then
$$
\|\theta_k-\theta^*\| \leq e^{-c_2 (k-k_0) \min_i r_{k,i}}\|\theta_{k_0}-\theta^*\|, \quad c_2:=\frac{\alpha \eta}{\sqrt{f(\theta_{k_0})+c})}.
$$
\end{theorem}
See Appendix A.3 for the proof.

\begin{remark} The above gradient estimate when using  $g(\theta)=\sqrt{f(\theta)+c}$
leads to 
\begin{equation*}
\min_{j<k}(\partial_i f(\theta_j))^2 
\leq \frac{1}{k}\sum_{j=0}^{k-1}(\partial_i f(\theta_j))^2
\leq \cfrac{2(f(\theta_0)+c)^{\frac{1}{2}}(F_k+c)}{\eta k r_{k,i}},\quad F_k=max_{j<k}f(\theta_j).
\end{equation*}
\end{remark}
As argued above, our convergence rates in Theorem \ref{thm3} depend also on the behavior of $r_{k, i}$, about which  we have the following result. 
\begin{lemma}\label{tau2}
Under the same assumptions as in Lemma \ref{tau}, for AEGD (\ref{ee1}) we have for $i\in [n]$, 
$$
r_{k,i}> r_i^* > 
g(\theta^*)(1-\eta/\tilde\tau),
\quad
\tilde\tau:=\frac{2g(\theta^*)}{nL_g (g(\theta_0))^2}.
$$
This implies that $r_i^*>0$ if $\eta <\tilde \tau$. 
\end{lemma}
We include a proof in Appendix A.2.

\section{Stochastic AEGD}
This section presents a stochastic AEGD algorithm for the unconstrained finite-sum optimization problem: 
\begin{equation}\label{dloss}
\min_{\theta} \left\{ f(\theta)=\frac{1}{m} \sum_{i=1}^m f_i(\theta)\right\}.    
\end{equation}
We use $\theta^*$ to denote a minimizer of $f(\theta)$, and assume that $f_i$ is differentiable and lower bounded so that 
$f_i(\theta)>-c$, for $i\in[m]$,
for some $c>0$. This problem is prevalent in machine learning tasks where $\theta$ corresponds to the model parameters, $f_i(\theta)$ represents the loss on the training point $i$ and the aim is to minimize the average loss $f(\theta)$ across points. When $m$ is large, SGD or its variants are preferred for solving (\ref{dloss}) mainly because of their cheap per iteration cost. To present a stochastic algorithm, we 
use $v_k$ to denote a random search direction at $k$-th step.   
Here we give only the element-wise version of our stochastic AEGD algorithm. 
\begin{algorithm}
\caption{Stochastic AEGD. Good default setting for parameters are $c=1$ and $\eta=0.1$. }
\label{alg:AEGD}
\begin{algorithmic}[1] 
\Require  $\{f_i({\theta})\}_{i=1}^m$, $\eta$: the step size, $\theta_0$: initial guess of $\theta$, and $T$: the total number of iterations.
\Require $c$: a parameter such that for any $i\in[m]$, $f_i({\theta})+c>0$ for all $\theta \in \mathbb{R}^n$, initial energy: $r_0=\sqrt{f_{i_0}(\theta_0)+c}{\bf 1}$
\For{$k=0$ to $k-1$}
\State $v_k:=\nabla f_{i_k}(\theta_k)/\big(2\sqrt{f_{i_k}(\theta_k)+c}\big)$ 
($i_k$ is a random sample from $[m]$ at step $k$) %
\State $r_{k+1} = r_k/(1+2\eta v_k\odot v_k)$ (update energy)
\State $\theta_{k+1} = \theta_k - 2\eta r_{k+1}\odot v_k$
\EndFor
\State \textbf{return} $\theta_k$
\end{algorithmic}
\end{algorithm}
\begin{remark}
If at each iteration step, a mini-batch of training data were selected, Algorithm \ref{alg:AEGD} still applies if $f_{i_k}(\theta_k)$ is replaced by
$
\frac{1}{b}\sum_{i \in B_k} f_i(\theta_k),
$
where $B_k$ denotes a randomly selected subset of $[m]$ of size $b$ at step $k$ with $b\ll m$. 
\end{remark}
\red{
}
To allow for any form of minibatching we use the arbitrary sampling notation 
$$
f_\xi(\theta)=\frac{1}{m}\sum_{j=1}^m \xi_j f_j(\theta),
$$
where $\xi=(\xi^1, \cdots, \xi^m) \in \mathbb{R}^m_+$ is a random sampling vector (drawn from some distribution) such that $\mathbb{E}[\xi^j]=1$ for $j\in[m].$ The element-wise update rule for the stochastic AEGD with arbitrary sampling can be reformulated as
\begin{subequations}\label{ee1a}
\begin{align}
& v_{k}=\nabla f_{\xi_k}(\theta_k)/(2\sqrt{f_{\xi_k}(\theta_k)+c}),\\
& r_{k+1,i}-r_{k,i}  = v_{k,i}(\theta_{k+1,i}-\theta_{k,i}), \quad r_{0,i} =\sqrt{f_{\xi_0}(\theta_0)+c}\\
& \theta_{k+1,i}=\theta_{k,i}-2\eta r_{k+1,i}v_{k,i} \quad k=0,1,2,\cdots.
\end{align}
\end{subequations}
It follows immediately from the definition of sampling vector $\xi$ that 
$$
\mathbb{E}[f_\xi(\theta)]=f(\theta), \quad \mathbb{E}[\nabla f_\xi(\theta)]=\nabla f(\theta),
$$
which means that we still have access to unbiased estimates of $f$ and its gradient. Of particular interest is the minibatch sampling: $\xi\in \mathbb{R}^m_+$ is a b-minibatch sampling if for every subset $B\in[m]$ with $|B|=b$  we have that
$$
\mathbb{P}\left[ 
\xi=\frac{m}{b}\sum_{i\in B}e_i
\right]=\frac{b!(m-b)!}{m!}.
$$
One can show by a double counting argument that if $\xi$ is a $b-$minibatch sampling, it is indeed a valid sampling with $\mathbb{E}[\xi^j]=1$ (see \cite{Get19}) 
and $\frac{1}{m}\sum_{j=1}^m \xi^j=1$. 

No matter how $v_k$ is defined, we have the following result. 
\begin{theorem}(Unconditional energy stability)
\label{thm1s1}  The stochastic AEGD of form (\ref{ee1a})
is unconditionally energy stable in the sense that for any step size $\eta>0$,
   \begin{equation}\label{srei+}
   \E[r^2_{k+1,i}]=\E[r^2_{k,i}] -\E[(r_{k+1,i} -r_{k,i})^2]- \eta^{-1} \E[(\theta_{k+1,i}-\theta_{k,i})^2],\quad i\in [n],
   \end{equation}
that is $\E[r_{k,i}]$ is strictly decreasing and convergent with $\E[r_{k,i}] \to r_i^*$ as $k\to \infty$, and also
\begin{equation}\label{srev1+}
   \lim_{k\to \infty}\E [(\theta_{k+1,i}-\theta_{k,i})^2]=0, \quad \sum_{j=0}^\infty\E[(\theta_{j+1,i}-\theta_{j,i})^2] \leq \eta (f(\theta_0)+c), \quad \forall i\in[n].
\end{equation}
\end{theorem}
The proof is entirely similar to that for Theorem \ref{thm1}, details are omitted.  

Due to the nonlinearity of $v_k$ in terms of the random variables, we only have 
$
 2\mathbb{E}\left[ v_k \sqrt{f_{\xi_k}(\theta_k)+c}\right]=\nabla f(\theta_k),
$
the estimate of convergence rates can be quite subtle.  Nevertheless, with mild assumptions  on $f_l$  we are able to establish the following.  
\begin{theorem}
\label{thm1s2} Suppose $\|\nabla f_l\|_\infty\leq G_\infty$, and $f_l+c \geq a>0$ for all $l\in [m]$. Then  stochastic AEGD of form (\ref{ee1a}) 
admits the following direction-wise estimate, for $i\in [n]$, 
$$
\frac{1}{k}\sum_{j=0}^{k}\E[(v_{j,i})^2] 
\leq \frac{C_i}{k\mathbb{E}[r_{k,i}]},  \quad
C_i:= \frac{\mathbb{E}[r_{0,i}](2a+\eta G^2_\infty)}{4 a\eta}.
$$
\end{theorem}
 The proof is given in Appendix \ref{pf-saegd}. %

\section{Experimental Results}\label{ch:experiments}

We evaluate the deterministic AEGD \eqref{ee1} and stochastic AEGD (Algorithm \ref{alg:AEGD}) on several benchmarks for optimization, including convex and non-convex performance testing problems, k-means clustering, which has a non-smooth objective function, and convolutional neural networks on the standard CIFAR-10 and CIFAR-100 data sets. Overall, we show that AEGD is a versatile algorithm that can efficiently solve a variety of optimization problems. 

In all experiments, we fine tune the base step size for each algorithm to obtain its best performance. The base step size is given in respective plots, denoted as ``lr'', the momentum in GDM and stochastic GDM (SGDM) is set to $0.9$. For AEGD we found that the parameters that impacts performance the most were the base step size and step decay schedule in the deep learning experiments. We, therefore, leave the energy parameter $c$ with the default value $1$.

Our experiments show the following primary findings: (i) AEGD allows larger effective step size, hence converges much faster than GD; (ii) Empirically the performance of AEGD appears better than or at least comparable with (S)GDM: in full batch setting, AEGD typically displays rapid initial progress while GDM tends to overshoot and detour to the target; in stochastic setting, AEGD produces solutions that generalize better than SGDM when coupled with a learning rate decay schedule.


\subsection{Performance Testing Problems}
We begin with comparing AEGD with GD and GDM on two benchmark convex and non-convex performance testing problems for optimization algorithms. Consider searching the minima $x^*$  of the following quadratic function (a strongly convex function with $\alpha=2/100$ and $L=2$)
\begin{equation}\label{eq:quad}
f(x_1,x_2,...,x_{100})=\sum_{i=1}^{50}x^2_{2i-1}+\sum_{i=1}^{50}\frac{x^2_{2i}}{10^2},
\end{equation}
and the 2D Rosenbrock function (non-convex and $L$-smooth)
\begin{equation}\label{eq:rosenbrock}
  f(x_1,x_2)=(1-x_1)^2+100(x_2-x_1^2)^2.
\end{equation}

\begin{figure}[ht]
\begin{subfigure}[b]{0.33\linewidth}
\centering
\includegraphics[width=1\linewidth]{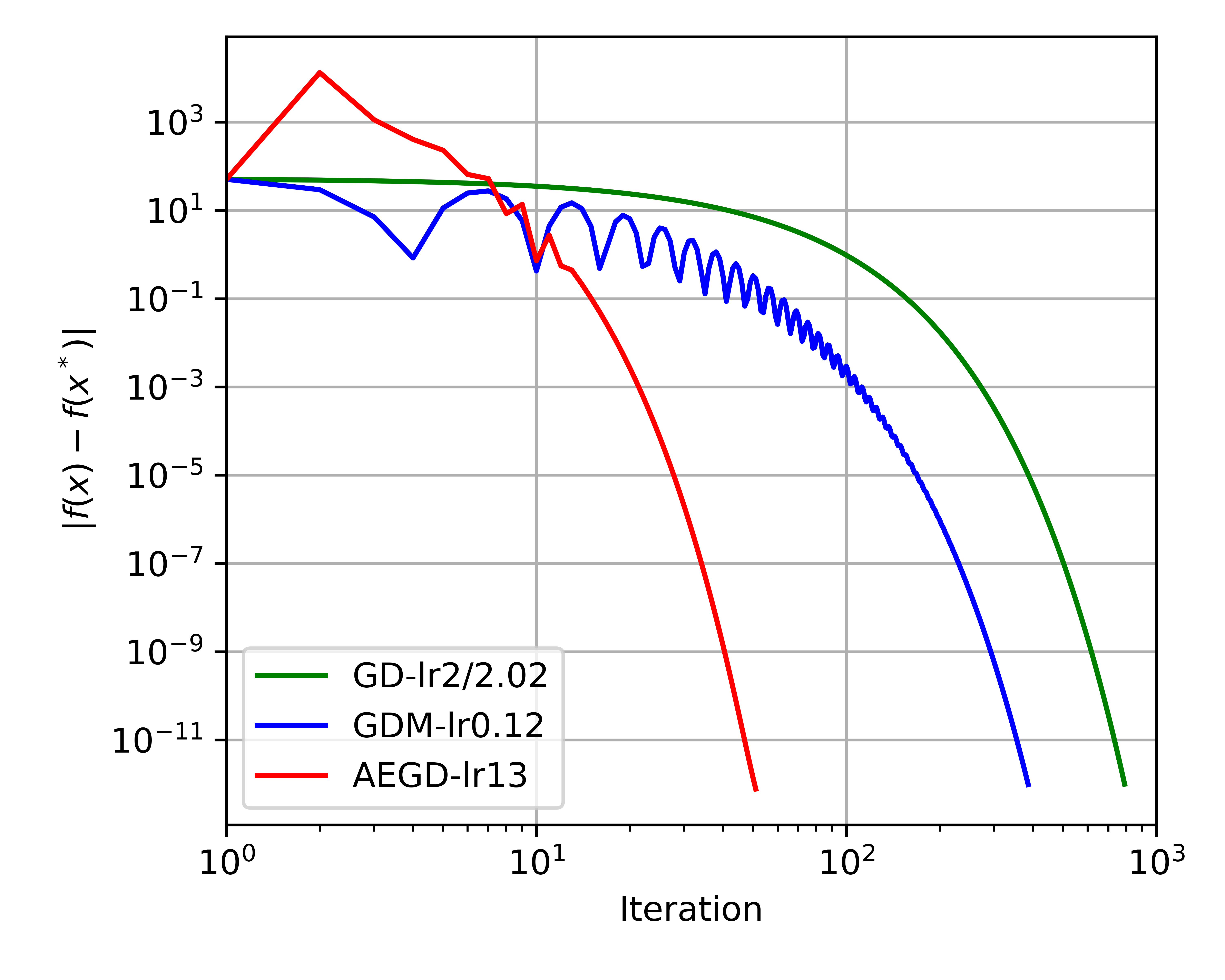}
\caption{}
\end{subfigure}%
\begin{subfigure}[b]{0.33\linewidth}
\centering
\includegraphics[width=1\linewidth]{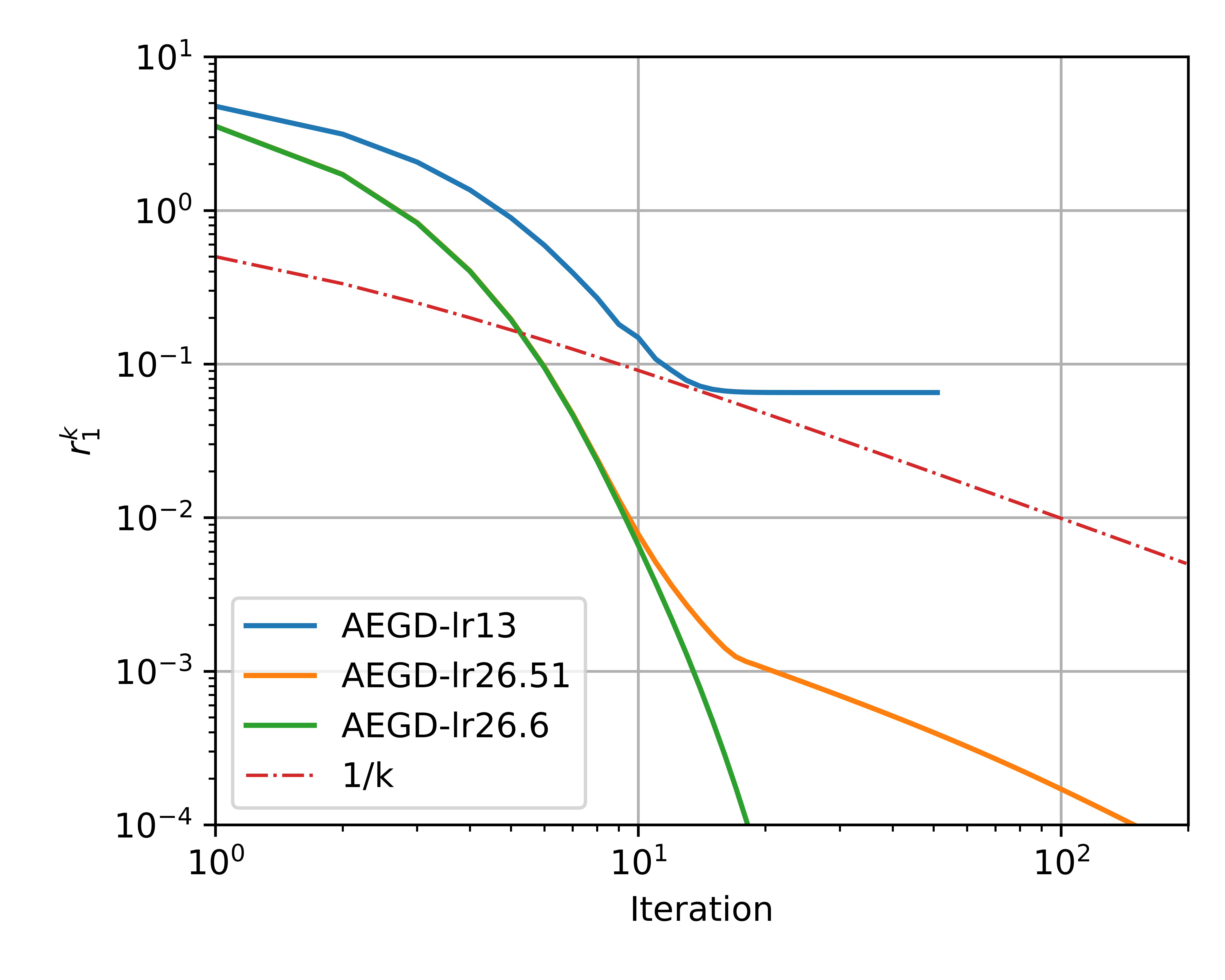}
\caption{}
\end{subfigure}%
\begin{subfigure}[b]{0.33\linewidth}
\centering
\includegraphics[width=1\linewidth]{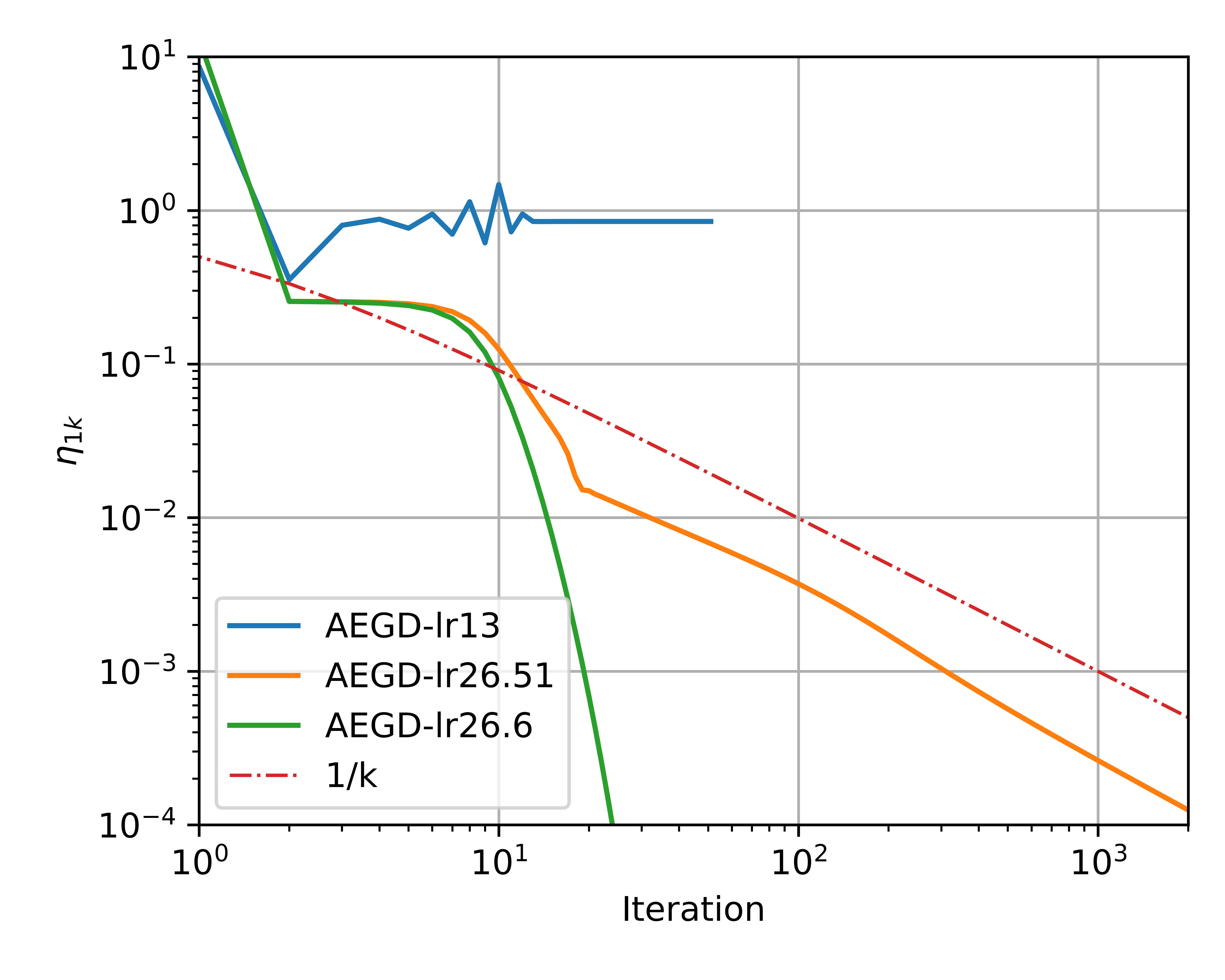}
\caption{} 
\end{subfigure}%
\captionsetup{format=hang}
\caption{The quadratic problem: (a) Optimality gap of different algorithms; (b) Behavior of $r_{k,1}$ for different $\eta$; (c) Behavior of $\eta_{1k}$ for different $\eta$.}
\label{fig:quad}
\end{figure}

\begin{figure}[ht]
\begin{subfigure}[b]{0.33\linewidth}
\centering
\includegraphics[width=1\linewidth]{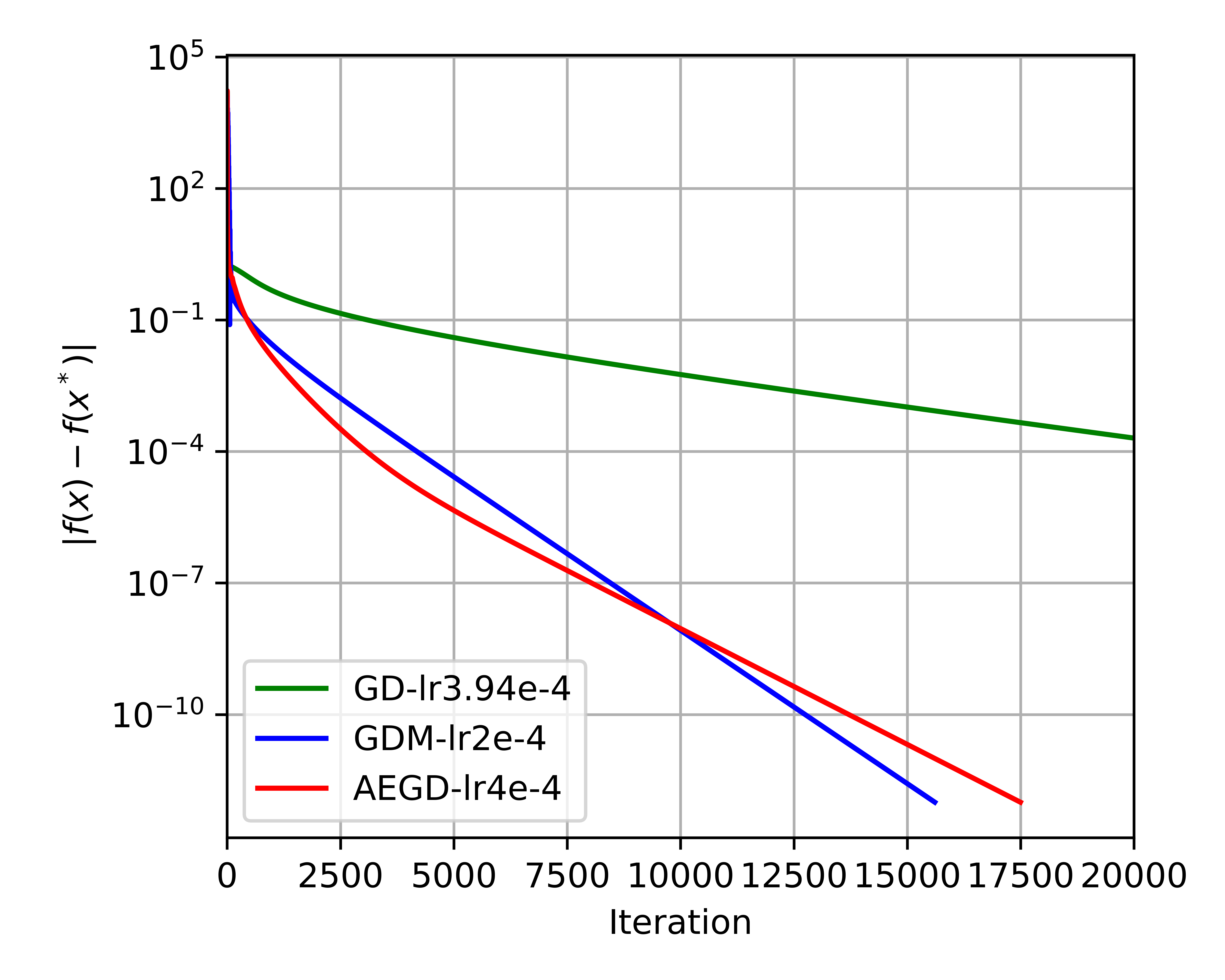}
\caption{}
\end{subfigure}%
\begin{subfigure}[b]{0.33\linewidth}
\centering
\includegraphics[width=1\linewidth]{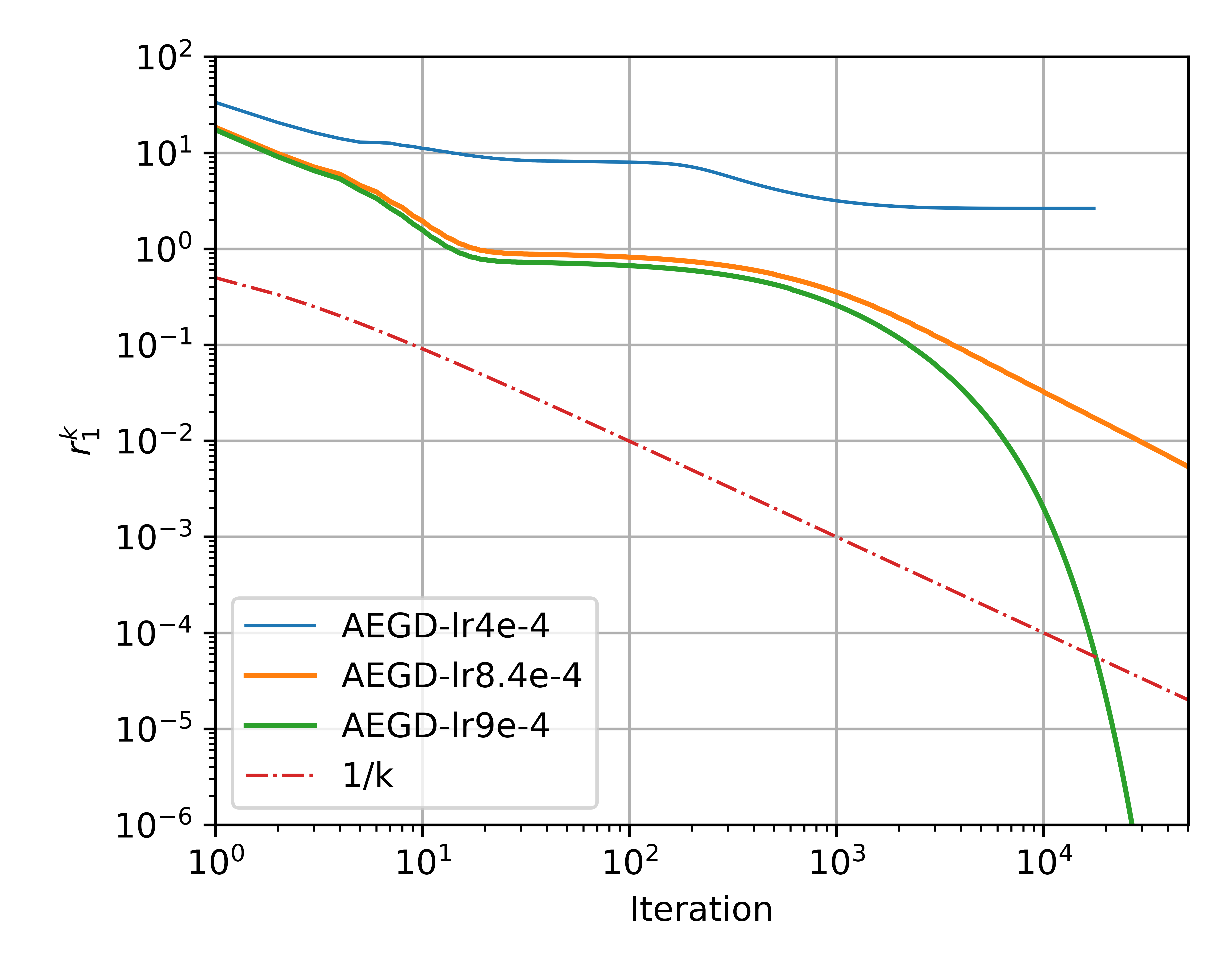}
\caption{}
\end{subfigure}%
\begin{subfigure}[b]{0.33\linewidth}
\centering
\includegraphics[width=1\linewidth]{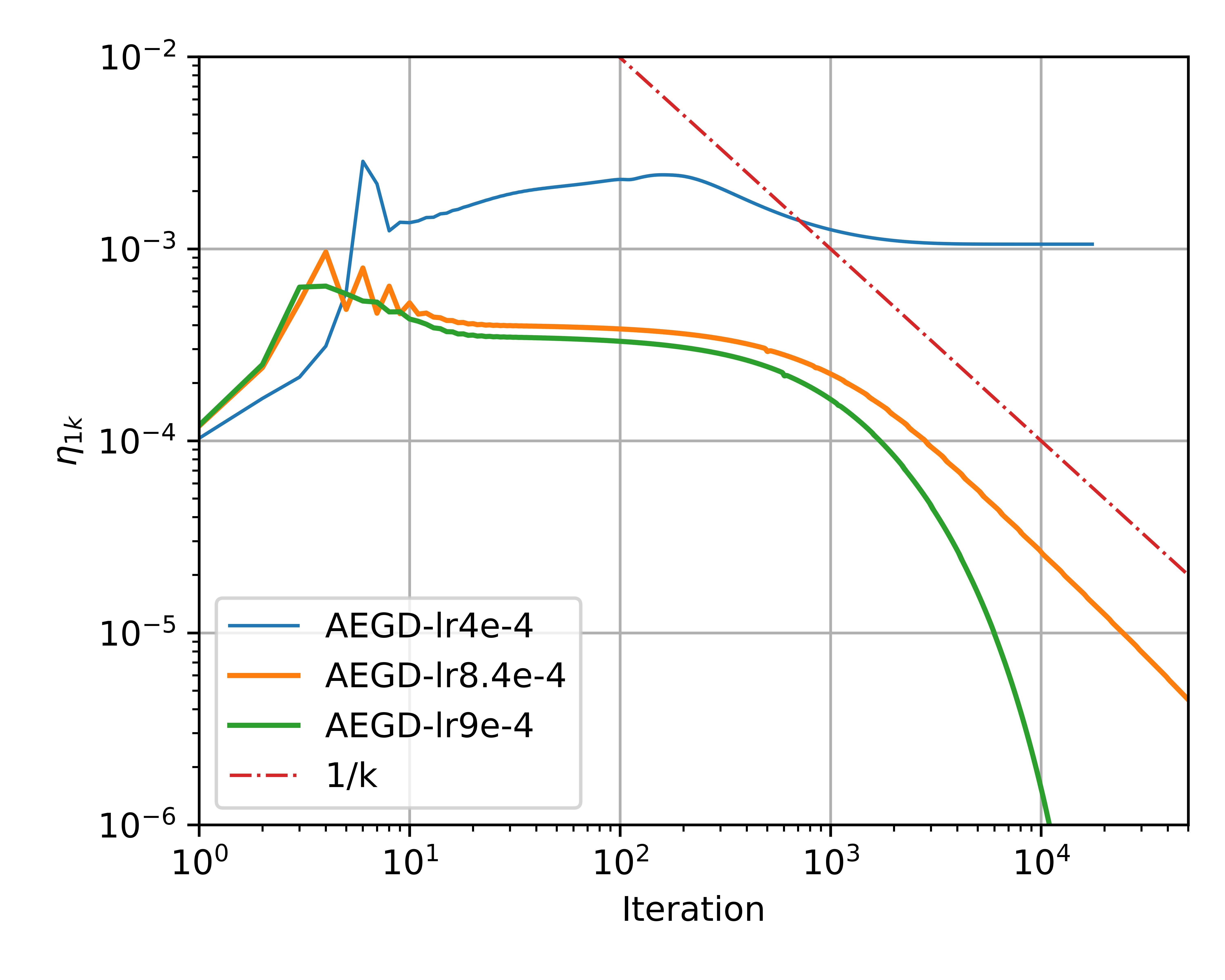}
\caption{}
\end{subfigure}%
\captionsetup{format=hang}
\caption{The Rosenbrock problem: (a) Optimality gap of different algorithms; (b) Behavior of $r_{k,1}$ for different $\eta$; (c) Behavior of $\eta_{1k}$ for different $\eta$.}
\label{fig:rosen}
\end{figure}

\begin{figure}[ht]
\centering
\includegraphics[width=0.4\linewidth]{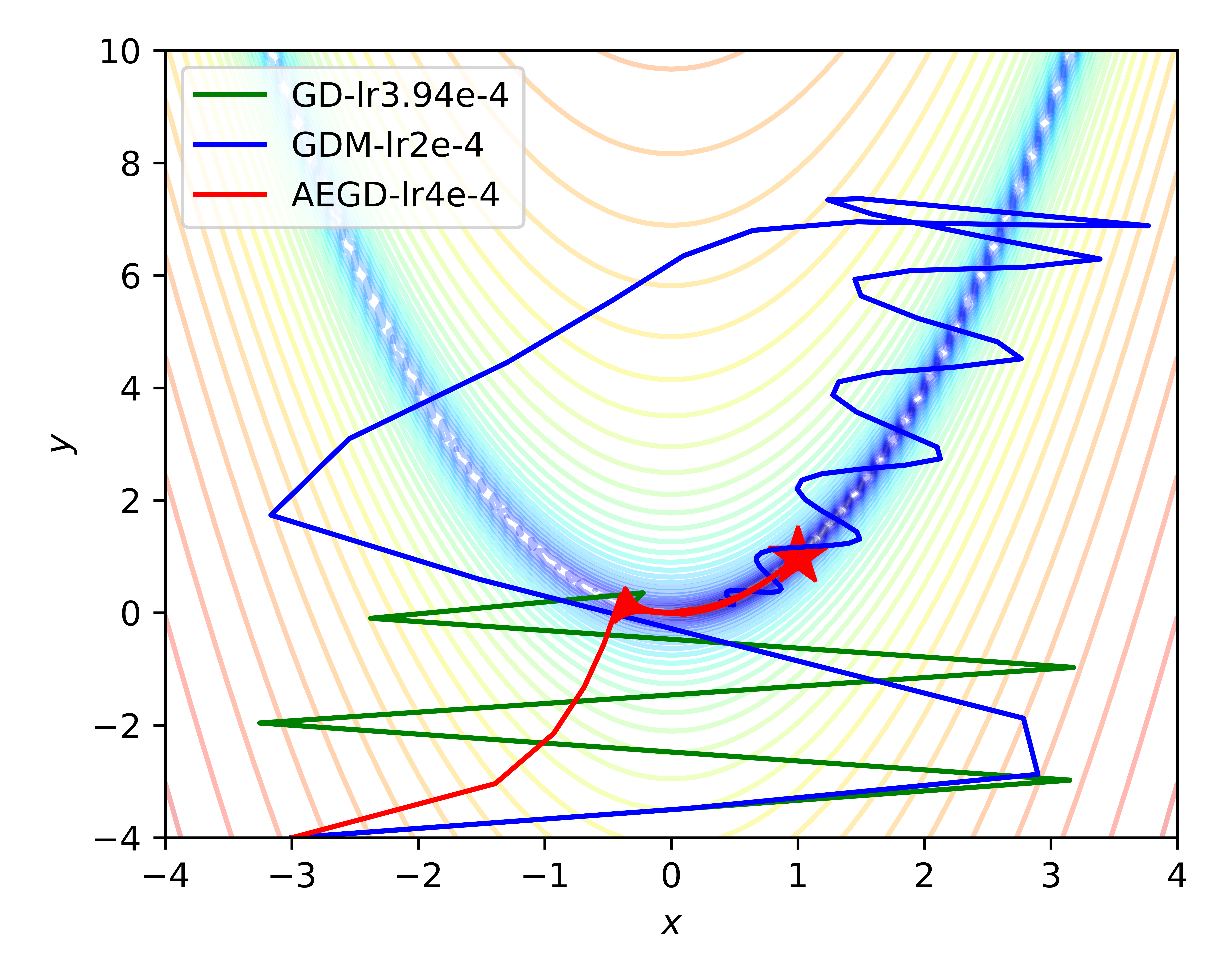}
\captionsetup{format=hang}
\caption{Trajectories of different optimization algorithms on the Rosenbrock function.}
\label{fig:rosen_paths}
\end{figure}

The initial points for (\ref{eq:quad}) and (\ref{eq:rosenbrock}) are set to $(1,1,...,1)$ and $(-3, -4)$, respectively. 
Numerical comparison results are in Figure \ref{fig:quad} (a) for the quadratic function and Figure \ref{fig:rosen} (a) for the Rosenbrock function. We see that AEGD converges much faster than both GD and GDM on the quadratic problem. For the Rosenbrock problem, the convergence of GD is very slow compared with GDM and AEGD. Though the number of iterations that AEGD needs to reach the minima is slightly more than GDM, AEGD makes a faster initial progress. This can be observed more clearly in Figure \ref{fig:rosen_paths}, which shows that GD and GDM tend to overshoot and detour to the minima while AEGD goes along a more direct path to the minima.

We also investigated how the basestep size $\eta$ affects the behavior of $r_k$ of AEGD on the two problems. (For both functions, $r_{k,1} =\min_{i} r_{k,i}$, hence we only show the behavior of $r_{k,1}$.) The results are presented in Figure \ref{fig:quad} (b) for the quadratic problem and Figure \ref{fig:rosen} (b) for the Rosenbrock problem. From these results, there appears to exist a threshold $\tilde\eta$ such that when $\eta<\tilde\eta$, $r^*>0$; when $\eta>\tilde\eta$, $r^*=0$; when $\eta=\tilde\eta$, $r_k$ converges to $0$ at the speed rate of $O(1/k)$; and when $\eta\leq\tilde\eta$, AEGD converges to the minima. We conjecture that such critical threshold phenomenon for the basestep size $\eta$ should hold true for more general objective functions.

Specifically, for the quadratic problem: $\tilde\eta \sim 26.51$; for the Rosenbrock problem, $\tilde\eta \sim 8.4e-4$. Based on our tests, the largest step size GD allows to ensure convergence is $1$ for the quadratic problem and $\sim3.94e-4$ for the Rosenbrock problem. In both cases, $\tilde\eta$ is much larger than the admissible step sizes for GD.


\subsection{K-means Clustering}

We consider $k$-means clustering problem for a set of data points $\{\bm{p}_i\}^m_{i=1}$ in $\mathbb{R}^d$ with $K$ centroids $\{\bm{x}_j\}^K_{j=1}$. Denote $\bm{x}=[\bm{x}_1, \cdots, \bm{x}_K]\in \mathbb{R}^{Kd}$,  we seek to minimize the quantization error:
\begin{equation}\label{knn_loss}
\min_{\bm{x}\in \mathbb{R}^{Kd}}\left\{f(\bm{x}): =\frac{1}{2m}\sum_{i=1}^m \min_{1\le j\le K} \|\bm{x}_j-\bm{p}_i\|^2 \right\}.
\end{equation}
If a data point $\bm{p}_i$ has more than one distinct nearest centroids, we assign $\bm{p}_i$ to one of them randomly. We define the gradient of $f$ (presented in Appendix \ref{ap:kmeans}) as \cite{BB95} did when applying the gradient-based method. In this example,  the Iris data set, which contains $150$ four-dimensional data samples from $3$ categories, is used to compare the robustness to initialization of three methods: GD, AEGD, and expectation maximization (EM) \cite{L82}.

\begin{figure}[ht]
\begin{subfigure}[b]{0.25\linewidth}
\centering
\includegraphics[width=1\linewidth]{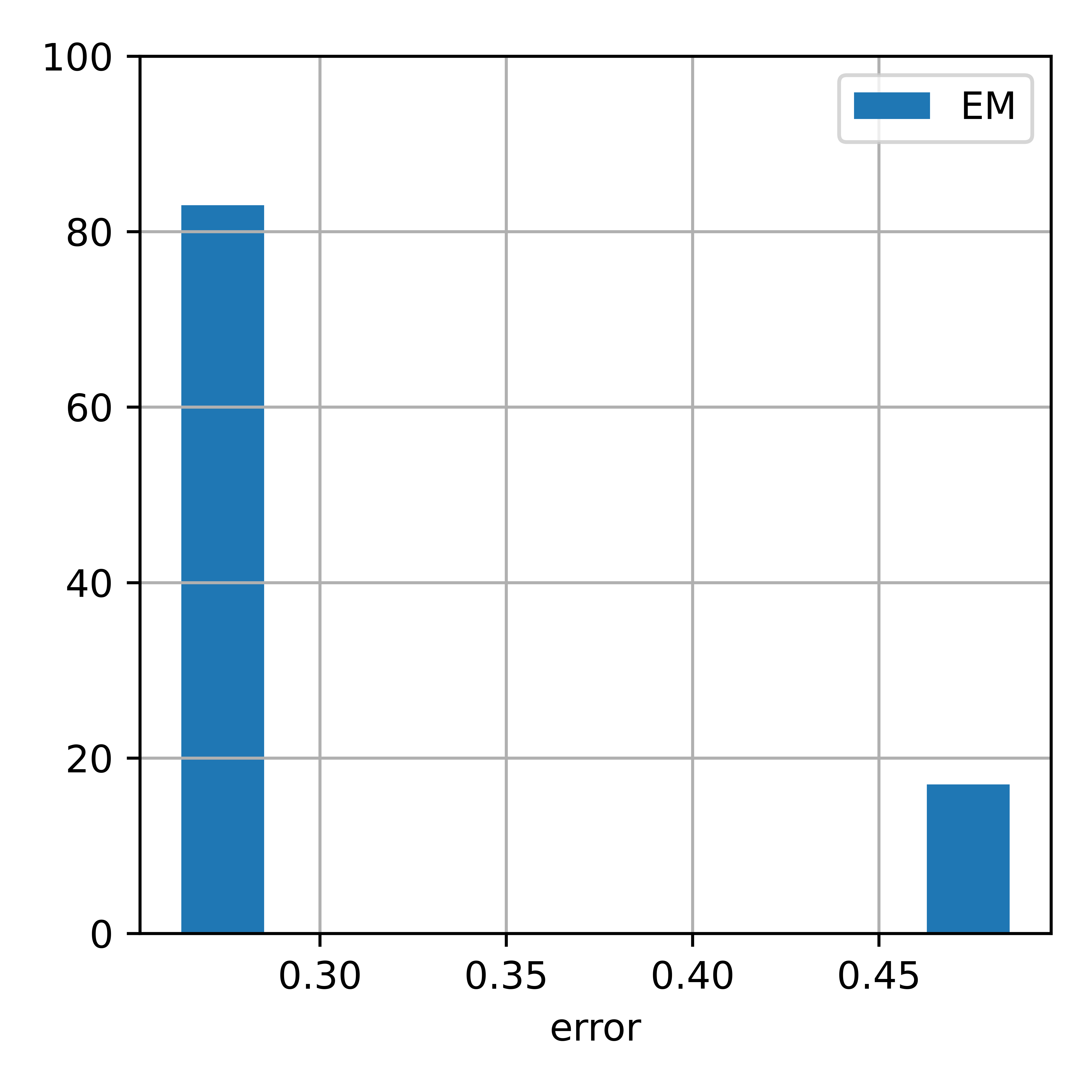}
\caption{}
\end{subfigure}%
\begin{subfigure}[b]{0.25\linewidth}
\centering
\includegraphics[width=1\linewidth]{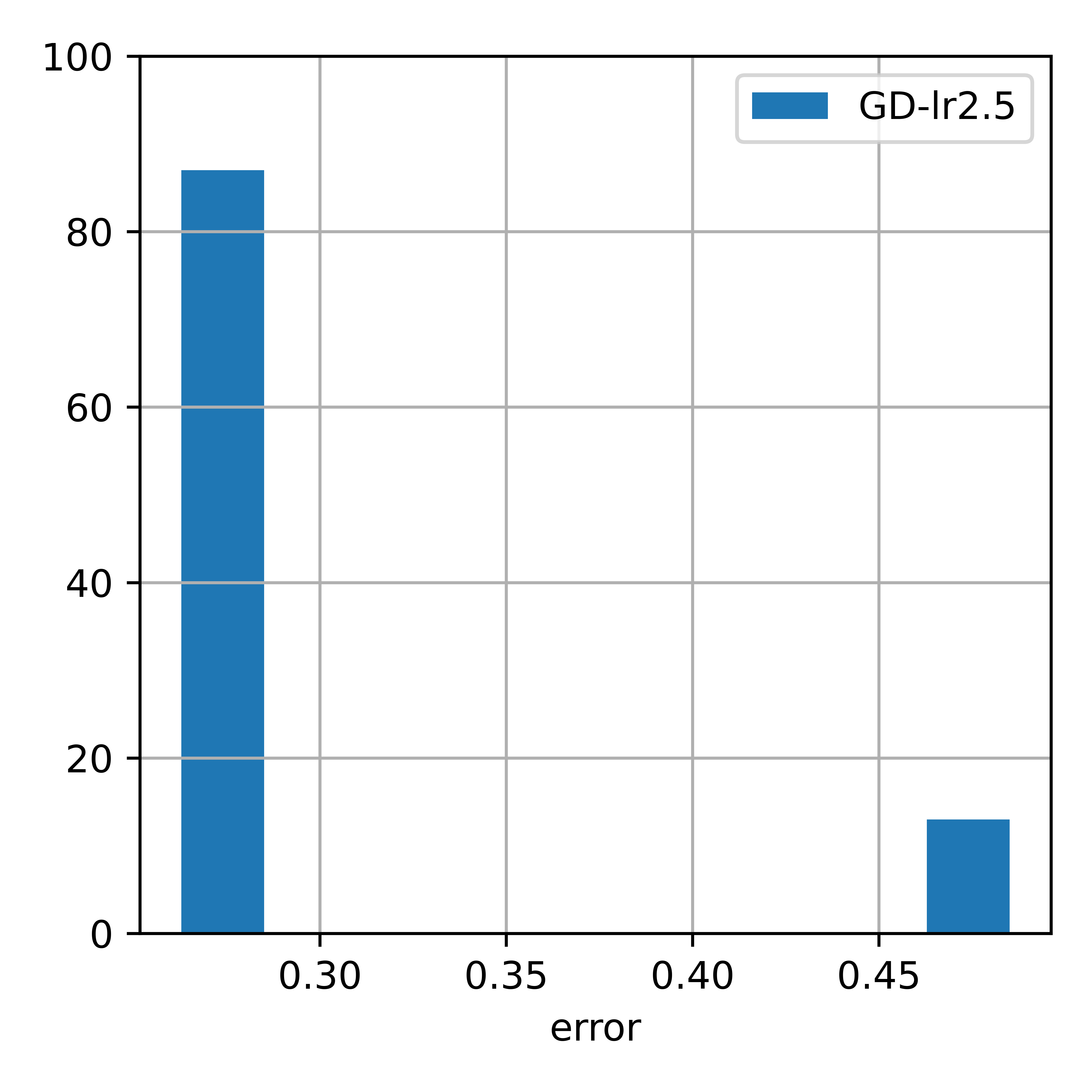}
\caption{}
\end{subfigure}%
\begin{subfigure}[b]{0.25\linewidth}
\centering
\includegraphics[width=1\linewidth]{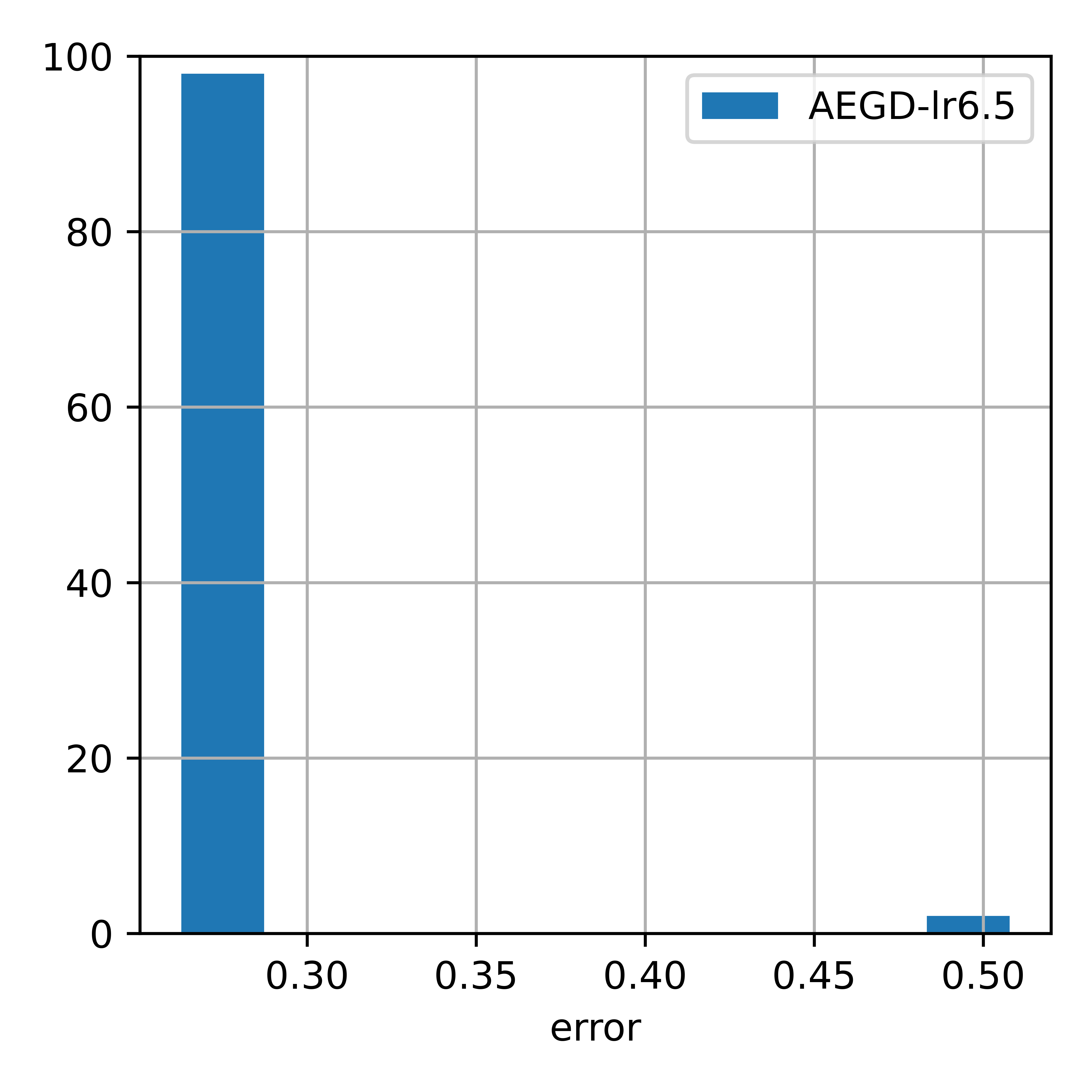}
\caption{}
\end{subfigure}%
\begin{subfigure}[b]{0.25\linewidth}
\centering
\includegraphics[width=1\linewidth]{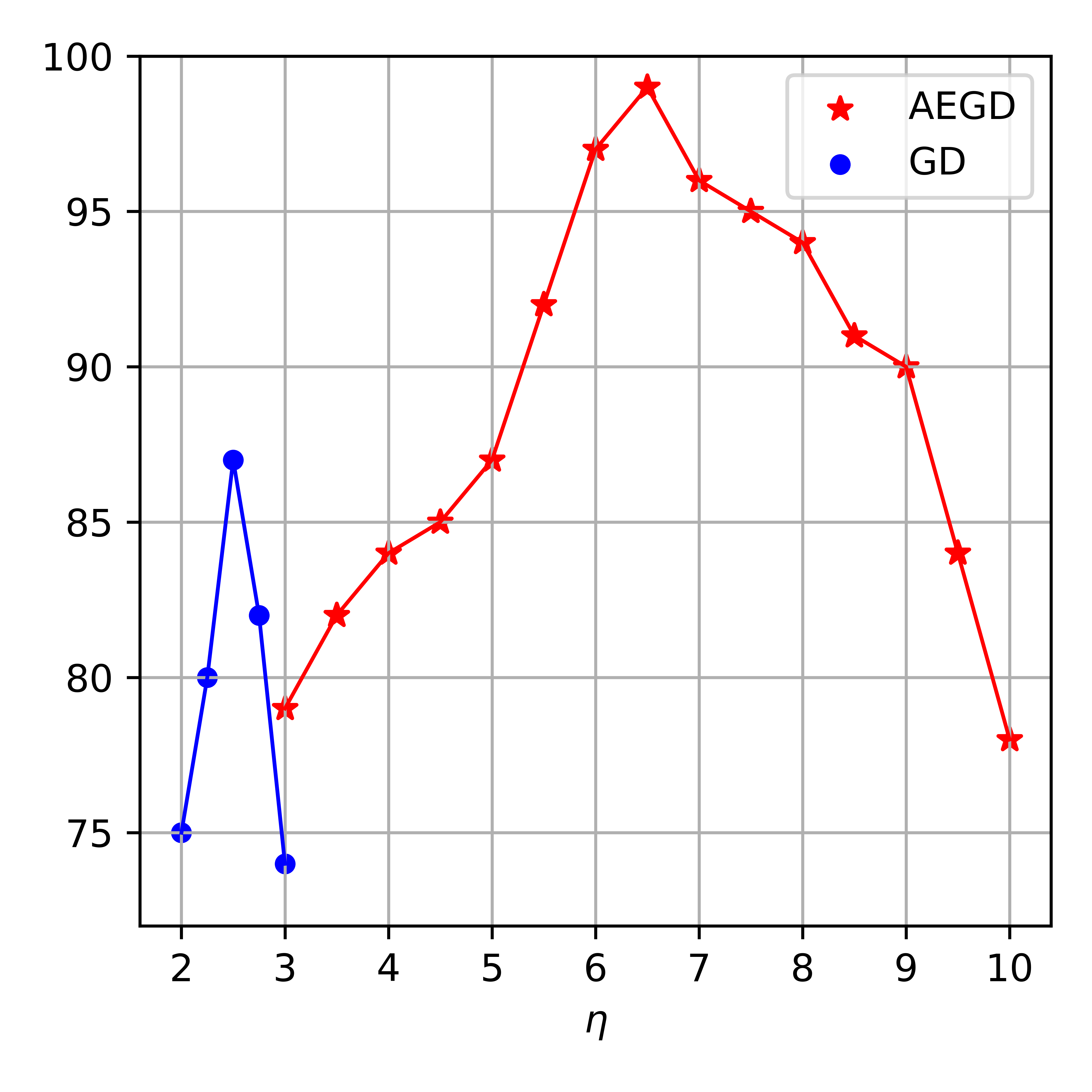}
\caption{}
\end{subfigure}%
\captionsetup{format=hang}
\caption{The histogram of the quantization error of $k$-means on Iris trained by EM (a), GD (b) and AEGD (c) over $100$ independent experiments. (d) Frequency of GD and AEGD achieve the improved minimum valued at $\sim 0.26$ in $100$ runs with different base step sizes.}
\label{fig:kmeans}
\end{figure}

Figure \ref{fig:kmeans} (a) (b) (c) present the frequency of error given by the three methods in $100$ runs. In each run, the initial centroids are selected from the data set randomly. We see that though there are chances for all the three methods to get stuck at a local minimum whose value is $\sim 0.48$, AEGD managed to locate an improved minimum valued at $\sim 0.26$ with the highest probability.

We also present the frequency of GD and AEGD achieving the improved minimum in $100$ runs with different base step sizes in Figure \ref{fig:kmeans} (d). We see that compared with GD, AEGD allows a larger set of base step size to achieve the improved minimum with much higher probability, with $\eta\sim 6.5$ being the optimal choice.

\subsection{Convolutional Neural Networks}
First, we should point out that the generalization capability of AEGD in training deep neural networks remains to be further understood. However, we would like to present some preliminary results to show the potential of the AEGD in this aspect.

\begin{figure}[ht]
\begin{subfigure}[b]{0.5\linewidth}
\centering
\includegraphics[width=1\linewidth]{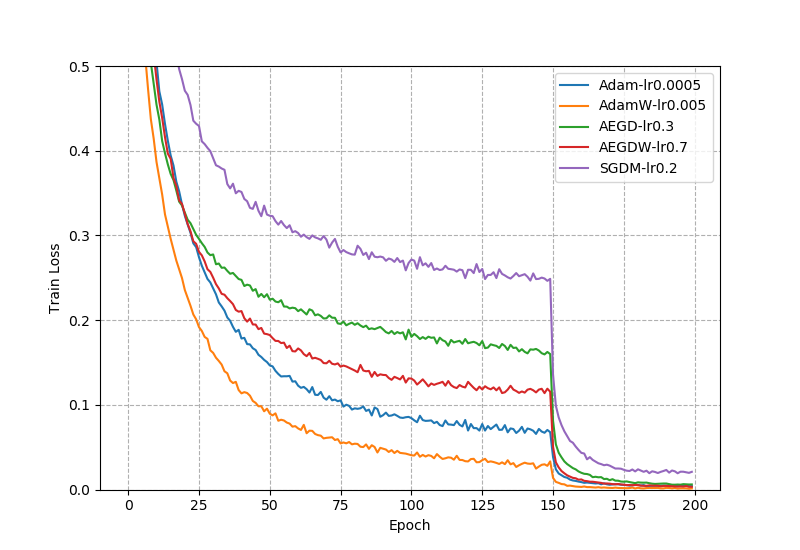}
\caption{Training loss, ResNet-56, CIFAR-10}
\end{subfigure}%
\begin{subfigure}[b]{0.5\linewidth}
\centering
\includegraphics[width=1\linewidth]{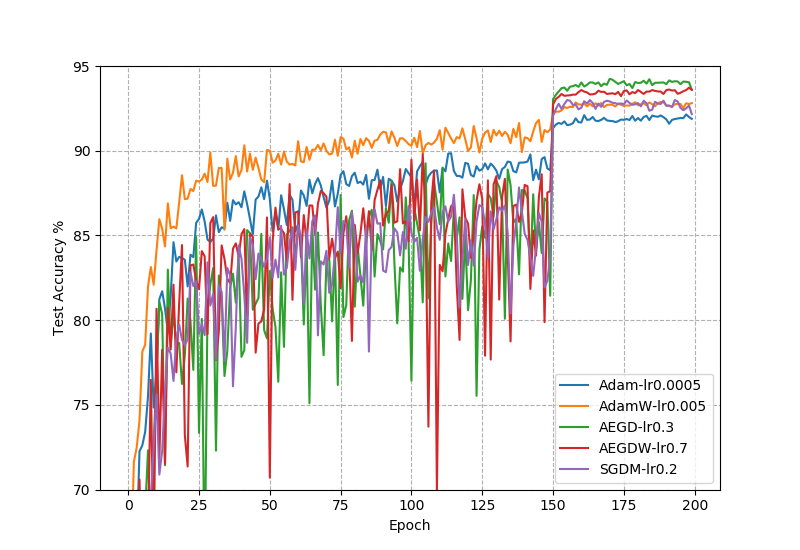}
\caption{Test accuracy, ResNet-56, CIFAR-10}
\end{subfigure}%
\newline
\begin{subfigure}[b]{0.5\linewidth}
\centering
\includegraphics[width=1\linewidth]{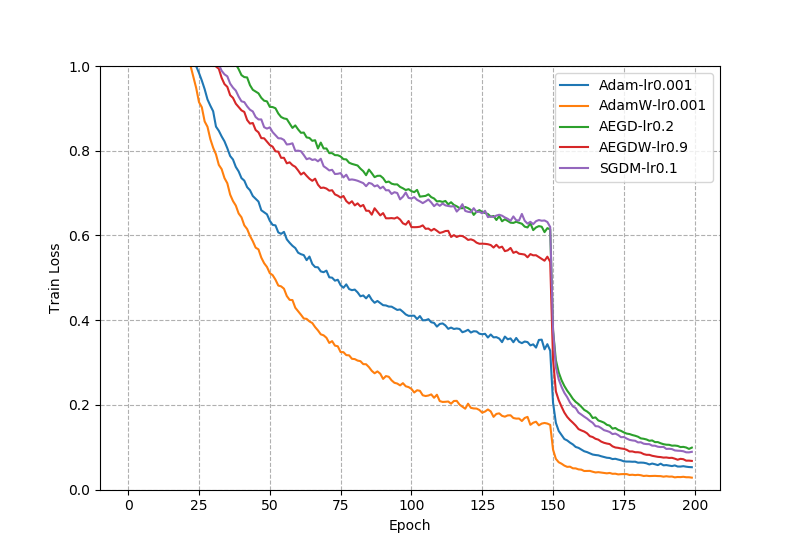}
\caption{Training loss, SqueezeNet, CIFAR-100}
\end{subfigure}%
\begin{subfigure}[b]{0.5\linewidth}
\centering
\includegraphics[width=1\linewidth]{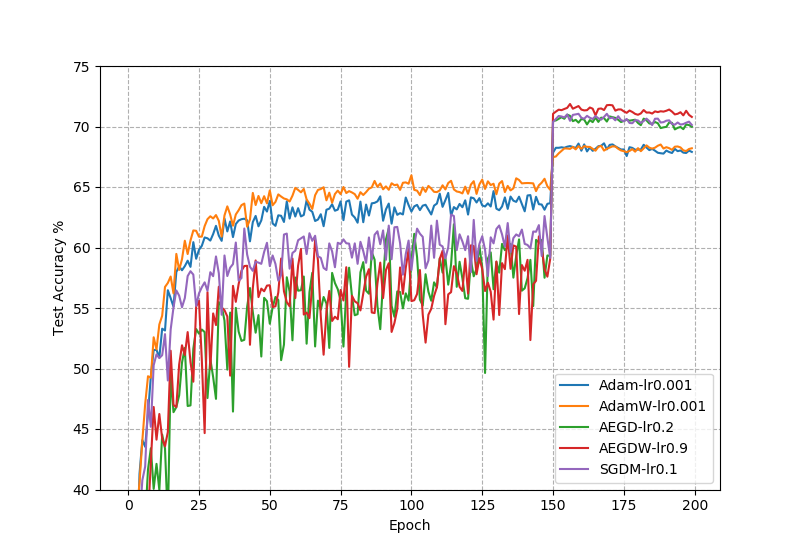}
\caption{Test accuracy, SqueezeNet, CIFAR-100}
\end{subfigure}%
\captionsetup{format=hang}
\caption{Training loss and test accuracy for ResNet-56 on CIFAR-10 and SqueezeNet on CIFAR-100}
\label{fig:cifar}
\end{figure}

Using ResNet-56 \cite{HZ+16} and SqueezeNet \cite{IH+16}, we consider the task of image classification on the standard CIFAR-10 and CIFAR-100 datasets. In our experiments, we employ the fixed budget of 200 epochs and reduce the learning rates by 10 after 150 epochs, with a minibatch size of $128$ and weight decay of $1\times10^{-4}$. We recall that weight decay has been a standard trick in training neural networks \cite{BC96, KH92}. For SGD, it can be interpreted as a form of $L^2$ regularization, but for adaptive algorithms such as Adam, careful implementation techniques are often needed \cite{LH19, WR17, ZWXG18}.  Inspired by \cite{LH19}, we introduce an AEGD-specific weight decay algorithm, called AEGDW. The algorithm for AEGDW and the initial set of step sizes are presented in Appendix \ref{ap:cifar}.

Our experimental results, as given in Figure \ref{fig:cifar}, show that Adam and AdamW perform better than other algorithms early. But by epoch $150$ when the learning rates are decayed, AEGD and AEGDW significantly outperform other methods in generalization. In the above two experiments, AEGD(W) even surpasses SGDM by $1\%$ in test accuracy. 

\section{Conclusions}\label{ch:conclusions}
Inspired by the IEQ approach for gradient flows in the form of time-dependent partial differential equations, we proposed AEGD in both global and element-wise form, as a new algorithm for first-order gradient-based optimization of non-convex objective functions. This simple update with an auxiliary energy variable is easy to implement, proven to be unconditionally energy stable irrespective of the step size, and features energy-dependent convergence rates with mild conditions on the step size. It is suitable for general objective functions, as long as they are bounded from below.  Numerical examples ranging from performance test problems, K-means clustering, to CNNs, all demonstrate the advantages of the proposed AEGD: it enjoys rapid initial progress and faster convergence than GD, is robust with respect to initial data, and generalizes better for deep learning problems. Overall, the results presented in this paper suggest that further research of AEGD could prove useful.

\appendix

\section{Technical Proofs}
In this section, we prove Theorem \ref{thm2} and Theorem \ref{thm3} from Section \ref{con}.
\subsection{Proof of Theorem \ref{thm2}}
We rewrite AEGD (\ref{e1+}) as
\begin{equation}\label{ee+}
 \theta_{k+1}= \theta_k- \eta_k \nabla f(\theta_k),\quad \eta_k:=\eta\frac{r_{k+1}}{g(\theta_k)}.
 \end{equation}
This is in the form of the usual GD with a variable step size $\eta_k$.\\
(i) 
Using scheme \eqref{e1} we have
$$
r_{j+1}-r_j = \nabla g(\theta_j) \cdot (\theta_{j+1}-\theta_j) = -2\eta r_{j+1}\|\nabla g(\theta_j)\|^2.
$$
Take a summation over $j$ from $0$ to $k-1$ gives
$$
r_0 - r_k = 2\eta\sum_{j=0}^{k-1}r_{j+1}\|\nabla g(\theta_j)\|^2 
\geq 2\eta r_k\sum_{j=0}^{k-1}\|\nabla g(\theta_j)\|^2
$$
Use $r_0=g(\theta_0)$ and $r_k$ is strictly decreasing to get
$$
k \min_{j<k}\|\nabla g(\theta_j)\|^2 
\leq \sum_{j=0}^{k-1}\|\nabla g(\theta_j)\|^2
\leq \frac{g(\theta_0)}{2\eta r_k}.
$$
Next we prove (i), (ii) and (iii), taking $k_0=0$, for simplicity in presentation.\\
(i) By L-smoothness of $f$ and scheme (\ref{ee+}), we have
\begin{equation}\label{df}
\begin{aligned}
f(\theta_{k+1}) -f(\theta_k) 
& \leq \nabla f(\theta_k)\cdot (\theta_{k+1}-\theta_k)+\frac{L}{2}\|\theta_{k+1}-\theta_k\|^2\\
& \leq -\frac{\eta_k}{2}( 2- L\eta_k)\|\nabla f(\theta_k)\|^2 \leq -\frac{\eta_k}{2}\|\nabla f(\theta_k)\|^2
\end{aligned}    
\end{equation}
as long as $\eta_k\leq 1/L$. Here and in what follows, we denote $w_k=f(\theta_k)-f^*$, then the PL property reads 
$$
\frac{1}{2}\|\nabla f(\theta_k)\|^2\geq\mu(f(\theta_k)-f^*)=\mu w_k.
$$
With this property, \eqref{df} can be written as
$$
w_{k+1}-w_k\leq -\mu\eta_k w_k.
$$
This implies $w_{k+1}\leq (1-\mu\eta_k) w_k$. By induction, 
\begin{align*}
w_k &\leq 
\prod_{j=0}^{k-1}(1-\mu\eta_j)w_0
= \exp\left(\sum_
{j=0}^{k-1} {\rm log}(1-\mu\eta_j )\right)w_0  \leq \exp\left(-\mu \sum_
{j=0}^{k-1} \eta_j\right)w_0.
\end{align*}
Noticing that
$
\eta_j=\eta \frac{r_{j+1}}{g(\theta_j)}\geq \eta \frac{r_{k}}{g(\theta_0)},
$
we further get
$$
w_k \leq \exp(-c_0k r_k)w_0, \quad c_0=\frac{\mu \eta }{g(\theta_0)}. 
$$
To obtain the convergence of $\theta_k$, we use scheme \eqref{ee+} to rewrite \eqref{df} as
$$
f(\theta_k) - f(\theta_{k+1})\geq \frac{1}{2\eta_k}\|\theta_{k+1}-\theta_k\|^2\quad\Rightarrow\quad
w_k-w_{k+1}\geq \frac{1}{2\eta_k}\|\theta_{k+1}-\theta_k\|^2.  
$$
The PL property when combined with (\ref{ee+}) gives  
$$
\|\theta_{k+1}-\theta_k\|^2\geq 2\mu\eta^2_kw_k \quad\Rightarrow\quad \frac{1}{\sqrt{w_k}}\geq\frac{\sqrt{2\mu}\eta_k}{\|\theta_{k+1}-\theta_k\|}.
$$
Using the above two inequalities and noting that $w_k\geq w_{k+1}$, we have
\begin{align*}
\sqrt{w_k}-\sqrt{w_{k+1}}
&\geq \frac{1}{2\sqrt{w_k}}(w_k-w_{k+1})
\geq \frac{\sqrt{2\mu}}{4}\|\theta_{k+1}-\theta_k\|.
\end{align*}
Taking a summation over $k$ from $0$ to $\infty$ gives 
$$
\sum_{k=0}^\infty\|\theta_{k+1}-\theta_k\|\leq \frac{4}{\sqrt{2\mu}} \sqrt{w_0}.
$$
This yields (\ref{ctheta}), which ensures the convergence of $\{\theta_k\}$. 
\\
(ii) 
By the $L$-smoothness assumption, scheme (\ref{ee+})
and $\eta_k \leq  1/L$, as in (i) we obtain 
$$
f(\theta_{k+1}) \leq f(\theta_k)-\frac{1}{2\eta_k}\|\theta_{k+1}-\theta_k\|^2.
$$
Denoting $d_k:=\|\theta_k-\theta^*\|$, we proceed to obtain the convergence rate. By convexity of $f$,
$$
f(\theta_k) \leq f(\theta^*) +\nabla f(\theta_k)\cdot (\theta_k-\theta^*).
$$
These when combined lead to
\begin{align*}
w_{k+1} & \leq \frac{1}{2\eta_k}( 2\eta_k \nabla f(\theta_k)\cdot (\theta_k -\theta^*)-\|\theta_{k+1}-\theta_k\|^2)\\
& =\frac{1}{2\eta_k}( - 2(\theta_{k+1}-\theta_k) \cdot (\theta_k -\theta^*)-\|\theta_{k+1}-\theta_k\|^2)\\
&  =\frac{1}{2\eta_k} (d_k^2 -d_{k+1}^2).
\end{align*}
Upon summation over iteration steps, we have
$$
d_0^2-d_k^2=2\sum_{j=0}^{k-1} \eta_i w_{i+1} \geq w_k\sum_{j=0}^{k-1}\eta_i \geq 2k w_k \frac{\eta r_k}{g(\theta_0)}.
$$
Hence for $\max_{j<k} \eta_j \leq 1/L$, we have
$$
f(\theta_k)-f(\theta^*)=w_k \leq \frac{\|\theta_0-\theta^*\|^2 g(\theta_0)}{2k\eta r_k}.
$$
(iii) To obtain the convergence rate of $\theta_k$,  we proceed with
\begin{align*}
\theta_{k+1} -\theta^* & = \theta_k-\theta^* -\eta_k \nabla f(\theta_k) \\
& =\theta_k -\theta^* - \eta_k(\nabla f(\theta_k) -\nabla f(\theta^*))\\
&=\theta_k-\theta^* -\eta_k \left( \int_0^1 \nabla^2 f(\theta^*+s(\theta_k-\theta^*))ds=:A_k\right)
(\theta_k-\theta^*).
\end{align*}
Therefore,
\begin{align*}
d_{k+1}& =\|(I-\eta_k A_k)(\theta_k-\theta^*)\|\\
& \leq \| I-\eta_k A_k \| d_k  \leq \Pi_{i=0}^{k} \|I-\eta_i A_i\|d_0\\
& \leq \max_{\alpha \leq \lambda \leq L} \Pi_{i=0}^{k}|1-\eta_i \lambda|d_0.
\end{align*}
For  $\max_{j<k} \eta_j \leq \frac{2}{\alpha +L}$, we have
\begin{align*}
d_{k}& \leq  \Pi_{i=0}^{k-1}(1-\eta_i \alpha )d_0
\leq e^{-c_2k r_k}d_0, \quad c_2=\frac{\alpha \eta }{g(\theta_0)}. 
\end{align*}
This completes the proof.
\subsection{Proof of Lemma \ref{tau2}}
The $L_g$-smoothness of $g$ implies that 
\begin{align*}
g(\theta_{j+1}) 
& \leq g(\theta_j)+\nabla g(\theta_j)\cdot (\theta_{j+1}-\theta_j) +\frac{L_g}{2}\|\theta_{j+1}-\theta_j\|^2 \\
& = g(\theta_j)+\sum_{i=1}^{n} \partial_i g(\theta_j)(\theta_{j+1,i}-\theta_{j,i}) + \frac{L_g}{2}\sum_{i=1}^{n}(\theta_{j+1,i}-\theta_{j,i})^2\\
& \leq  g(\theta_j) + \sum_{i=1}^{n}(r_{j+1,i}-r_{j,i}) + \frac{\eta L_g}{2}\sum_{i=1}^{n}(r^2_{j,i}-r^2_{j+1,i}).
\end{align*}
Take summation over $j$ from $0$ to $k-1$ and use $r_{0,i}=g(\theta_0)$,
so that 
\begin{align*}
g(\theta_k)-g(\theta_0) 
&\leq \sum_{i=1}^{n}r_{k,i}-ng(\theta_0) + \frac{\eta L_g }{2} \left( 
n(g(\theta_0))^2-\sum_{i=1}^{n}r^2_{k,i} 
\right).  
\end{align*}
Using 
$g(\theta_k)\geq g(\theta^*)$, we have for any $k$, 
$$
-\sum_{i=1}^{n}r_{k,i} +g(\theta^*) + (n-1)g(\theta_0) - \frac{n \eta L_g}{2} (g(\theta_0))^2 \leq 0. 
$$
Passing to the limit as $k \to \infty$, we also have
$$
-\sum_{i=1}^{n}r_i^* + g(\theta^*) + (n-1)g(\theta_0) - \frac{n \eta L_g}{2} (g(\theta_0))^2 \leq 0. 
$$
From this we get 
$$
- \left(\min_i r_i^*+(n-1)g(\theta_0)\right) + g(\theta^*) + (n-1)g(\theta_0) - \frac{n \eta L_g}{2} (g(\theta_0))^2 \leq 0,
$$
which can be reduced to
$$
- \min_i r_i^* + g(\theta^*)(1-\eta/\tilde\tau) \leq 0,\quad \tilde\tau:=\frac{2g(\theta^*)}{nL_g (g(\theta_0))^2}.
$$
Hence 
$$
\min_i r_i^* > g(\theta^*)(1-\eta/\tilde\tau).
$$ 
\subsection{Proof of Theorem \ref{thm3} }
We rewrite AEGD (\ref{ee1}) for $i\in [n]$ as
\begin{equation}\label{ee1+}
 \theta_{k+1,i}= \theta_{k,i}- \eta_{ik} \partial_i f(\theta_k),\quad \eta_{ik}:=\eta\frac{r_{k+1,i}}{g(\theta_k)}.
\end{equation}
(i) 
Using scheme \eqref{ri}, for $i\in[n]$ we have
$$
r_{j+1,i}-r_{j,i} = \partial_i g(\theta_j) (\theta_{j+1,i}-\theta_{j,i}) = -2\eta r_{j+1,i}(\partial_i g(\theta_j))^2.
$$
{
Take summation over $j$ from $0$ to $k-1$ gives
\begin{align*}
r_{0,i} - r_{k,i}
= 2\eta\sum_{j=0}^{k-1}r_{j+1,i}(\partial_i g(\theta_j))^2 
\geq 2\eta r_{k,i} \sum_{j=0}^{k-1}(\partial_i g(\theta_j))^2
\end{align*}
Using $r_{0,i}=g(\theta_0)$ and $r_{j,i}$ is strictly decreasing, we get
$$
k\min_{j<k}(\partial_i g(\theta_j))^2 
\leq \sum_{j=0}^{k-1}(\partial_i g(\theta_j))^2
\leq \frac{g(\theta_0)}{2\eta r_{k,i}}.
$$
}

Next we turn to prove (i), (ii) and (iii). For simplicity in presentation, we take $k_0=0$.\\
(i) By L-smoothness of $f$ and scheme (\ref{ee+}), we have
\begin{align*}
f(\theta_{k+1}) -f(\theta_k) 
& \leq \nabla f(\theta_k)\cdot (\theta_{k+1}-\theta_k)+\frac{L}{2}\|\theta_{k+1}-\theta_k\|^2\\
& \leq -\sum_{i=1}^{n}\frac{\eta_{ik}}{2}(\partial_i f(\theta_k))^2
\leq -\frac{\min_i\eta_{ik}}{2}\|\nabla f(\theta_k)\|^2
\end{align*}    
as long as $\eta_{ik}\leq 1/L$.
As in the proof for Theorem \ref{thm2} (i), we have
\begin{align*}
w_k
\leq \exp\left(-\mu \sum_
{j=0}^{k-1} \min_i\eta_{ij}\right)w_0.
\end{align*}
Noticing that
$
\min_i\eta_{ij}=\min_i\eta \frac{r_{j+1,i}}{g(\theta_j)}\geq \eta \frac{\min_i r_{k,i}}{g(\theta_0)},
$
we further get
$$
w_k \leq \exp(-c_0k \min_i r_{k,i})w_0, \quad c_0=\frac{\mu \eta }{g(\theta_0)}. 
$$
(ii) Using the L-smoothness assumption, scheme (\ref{ee1}) 
and $\eta_{ik} \leq  1/L$, as in (i) we obtain 
$$
f(\theta_{k+1}) \leq f(\theta_k)-\sum_{i=1}^n \frac{1}{2\eta_{ik}}(\theta_{k+1,i}-\theta_{k,i})^2.
$$
Denote $d_{ik}:=\theta_{k,i}-\theta_i^*$, then by convexity of $f$,
$$
f(\theta_k) \leq f(\theta^*) +\nabla f(\theta_k)\cdot (\theta_k-\theta^*).
$$
These when combined lead to
\begin{align*}
w_{k+1} & \leq \sum_{i=1}^n \frac{1}{2\eta_{ik}}
( 2\eta_{ik} \partial_i f(\theta_k) (\theta_{k,i} -\theta_i^*)-(\theta_{k+1,i}-\theta_{k,i})^2)\\
& =\sum_{i=1}^n \frac{1}{2\eta_{ik}}( - 2(\theta_{k+1,i}-\theta_{k,i}) (\theta_{k,i} -\theta_i^*)-(\theta_{k+1,i}-\theta_{k,i})^2)\\
&  =\sum_{i=1}^n \frac{1}{2\eta_{ik}} (d_{ik}^2 -d_{i,k+1}^2).
\end{align*}
That is
$$
\frac{2\eta}{g(\theta_k)} w_{k+1} \leq \sum_{i=1}^n \frac{1}{r_{k+1,i}} (d_{ik}^2 -d_{i,k+1}^2).
$$
This upon summation over iteration steps gives
\begin{align*}
\sum_{j=0}^{k-1}\frac{2\eta}{g(\theta_j)} w_{j+1} & \leq \sum_{j=0}^{k-1}\sum_{i=1}^n \frac{1}{r_{j+1,i}} (d_{ij}^2 -d_{i,j+1}^2) \\
& \leq \sum_{i=1}^n \left[ \frac{d^2_{i0}}{r_{1,i}} -\frac{d^2_{ik}}{r_{k,i}} +\sum_{j=0}^{k-1}
\left( \frac{1}{r_{j+1,i}}-\frac{1}{r_{j,i}}\right)d^2_{ij}\right]=:RHS.
\end{align*}
Since $f(\theta_j)$ is decreasing in $j$, so is $g(\theta_j)$. Hence
$$
\sum_{j=0}^{k-1}\frac{2\eta}{g(\theta_j)} w_{j+1}\geq  2k\eta \frac{w_k}{g(\theta_0)}.
$$
On the other hand, using  $r_{j,i} >r_{j+1,i}$, we have
\begin{align*}
RHS & \leq \sum_{i=1}^n \left[ \frac{d^2_{i0}}{r_{1,i}} -\frac{d^2_{ik}}{r_{k,i}} +\max_{j<k}d^2_{ij} \sum_{j=0}^{k-1} \left( \frac{1}{r_{j+1,i}}-\frac{1}{r_{j,i}}\right)\right] \\
& \leq \sum_{i=1}^n \left[ \frac{d^2_{i0}}{r_{1,i}} -\frac{d^2_{ik}}{r_{k,i}} +\max_{j<k}d^2_{ij}  \left( \frac{1}{r_{k,i}}-\frac{1}{r_{0,i}}\right)\right] \\
& \leq 2 \sum_{i=1}^n  \frac{\max_{j<k}d^2_{ij}}{r_{k,i}}.
\end{align*}
Hence
$$
w_k \leq \frac{g(\theta_0)}{k\eta}\sum_{i=1}^n  \frac{\max_{j<k}|\theta_{j,i} -\theta_i^*|^2}{r_{k,i}} \leq \frac{g(\theta_0)}{k\eta}  \frac{\max_{j<k}\|\theta_j -\theta^*\|^2}{ \min_i r_{k,i}}.
$$
(iii) To obtain the rate of $\theta_k$ converging to $\theta^*$,  we proceed with
\begin{align*}
\theta_{k+1,i} -\theta_i^* 
& =\theta_{k,i} -\theta_i^* - \eta_{ik}(\partial_i f(\theta_k) -\partial_i f(\theta^*))\\
&=\theta_{k,i}-\theta_i^* -\eta_{ik} \sum_{l=1}^n \left( \int_0^1 \partial_i\partial_l f(\theta^*+s(\theta_k-\theta^*))ds=:a_{il}^k\right)
(\theta_l^k-\theta_l^*)
\end{align*}
That is
$$
\theta_{k+1}-\theta^*=(I-{\rm diag}(\eta_k) A_k)(\theta_{k+1}-\theta^*),
$$
where ${\rm diag}(\eta_k)={\rm diag}(\eta_{1k},\cdots, \eta_{nk})$ and $A_k=(a_{il}^k)$. By matrix norm properties, we obtain
\begin{align*}
\|\theta_k-\theta^*\| =\Pi_{j=0}^{k-1} \left(\max_i \max_{\lambda \in [\alpha, L]} |1-\eta_{ik}\lambda| \right)\|\theta_0-\theta^*\|.
\end{align*}
Under the assumption $\eta_{ik}\leq \frac{2}{\alpha+L}$, we have
\begin{align*}
\|\theta_k-\theta^*\| & \leq  \Pi_{j=0}^{k-1}\max_i (1-\eta_{ij} \alpha )\|\theta_0-\theta^*\| \leq \exp\left(\sum_
{j=0}^{k-1} {\rm log}(1- \min_i \eta_{ij} \alpha)\right)\|\theta_0-\theta^*\| \\
& \leq \exp\left(-\alpha \sum_
{j=0}^{k-1} \min_i \eta_{ij}\right)\|\theta_0-\theta^*\|
 \leq e^{-k \alpha \eta \min_i r_{k,i}/ g(\theta_0)}\|\theta_0-\theta^*\|.
\end{align*}
In the last inequality we used the estimate
$\eta_{ij}=\eta \frac{r_{j+1,i}}{g(\theta_j)}
\geq \eta \frac{r_{k,i}}{g(\theta_0)}$ for $j<k.$
\subsection{Proof of Theorem \ref{thm1s2}}\label{pf-saegd}
{From
$
v_{j,i} = \partial_i f_{\xi}(\theta_j) /(2\sqrt{f_{\xi}(\theta_j)+c}),  
$
it follows  
$
(v_{j,i})^2 \leq 
\frac{G_\infty^2}{4a}. 
$
Using (\ref{ee1a}a, b), we obtain 
\begin{equation}\label{sri+}
r_{j,i} - r_{j+1,i} 
= - v_{j,i} (\theta_{j+1,i}-\theta_{j,i}) = 2\eta r_{j+1,i}(v_{j,i})^2= 
2\eta r_{j,i}\frac{(v_{j,i})^2}{1+2\eta(v_{j,i})^2},
\end{equation}
where we used 
$
r_{j+1,i}=\frac{r_{j,i}}{1+2\eta(v_{j,i})^2}. 
$
Taking expectation conditioned on $(\theta_j, r_j)$ in \eqref{sri+}, we have  
\begin{equation}\label{conep}
r_{j,i}- \mathbb{E}[r_{j+1,i}] 
= 2\eta r_{j,i}\mathbb{E}\bigg[\frac{(v_{j,i})^2}{1+2\eta(v_{j,i})^2}\bigg] 
\geq \frac{2\eta r_{j,i}}{1+2\eta G_\infty^2/(4a)}\mathbb{E}[(v_{j,i})^2].
\end{equation}
Rearranging and taking expectations to get
$$
\mathbb{E}[r_{j,i}]- \mathbb{E}[r_{j+1,i}]  \geq\frac{4a\eta }{2a +\eta G_\infty^2}\mathbb{E}[r_{j,i}]\mathbb{E}[(v_{j,i})^2].
$$
Summing over $j$ from $0$ to $k-1$ and using telescopic cancellation gives 
$$
\mathbb{E}[r_{0,i}]- \mathbb{E}[r_{k,i}]  
\geq\frac{4a\eta }{2a +\eta G_\infty^2}\sum_{j=0}^{k-1} \mathbb{E}[r_{j,i}]\mathbb{E}[(v_{j,i})^2]
\geq\frac{4a\eta }{2a +\eta G_\infty^2}\mathbb{E}[r_{k,i}]\sum_{j=0}^{k-1} \mathbb{E}[(v_{j,i})^2].  
$$
That is 
$$
k\min_{j<k}\mathbb{E}[(v_{j,i})^2] 
\leq \sum_{j=0}^{k-1} \mathbb{E}[(v_{j,i})^2]
\leq \frac{C_i}{\mathbb{E}[r_{k,i}]}, \quad
C_i:= \frac{\mathbb{E}[r_{0,i}](2a+\eta G_\infty^2)}{4 a\eta}.
$$
This completes the proof. 
}

\section{Additional Experimental Results and Implementation Details}
Here we provide additional experimental results and implementation details beyond those in Section \ref{ch:experiments}.

\subsection{K-means Clustering}\label{ap:kmeans}
By abuse of notation, we can define the `gradient' of $f$ at any point $\bm{x}$ as
\begin{equation}\label{knn_grad}
\nabla f(\bm{x})=\frac{1}{m}\left[\sum_{i\in \mathcal{C}_1}(\bm{x}_1-\bm{p}_i),...,\sum_{i\in \mathcal{C}_K}(\bm{x}_K-\bm{p}_i)\right]^{\top},
\end{equation}%
where $\mathcal{C}_j$ denotes the index set of the points that are assigned to the centroid $\bm{x}_j$.
With the definition of loss function in (\ref{knn_loss}) and its gradient so defined, we can apply gradient-based methods including AEGD to solve the $k-$means clustering problem.

\subsection{Convolutional Neural Networks}\label{ap:cifar}
The initial set of step sizes used for each algorithm are\\
SGDM: \{0.05, 0.1, 0.2, 0.3\},\\
Adam: \{1e-4, 3e-4, 5e-4, 1e-3, 2e-3\},\\
AdamW: \{5e-4, 1e-3, 3e-3, 5e-3\},\\
AEGD: \{0.1, 0.2, 0.3, 0.4\},\\
AEGDW: \{0.6, 0.7, 0.8, 0.9\}.

The results presented in Figure \ref{fig:cifar} show that Adam(W) does not generalize as well as SGDM and AEGD(W). Therefore, we only compare the generalization capability of AEGD(W) with SGDM for the remainder of the experiments.

\begin{figure}[ht]
\begin{subfigure}[b]{0.5\linewidth}
\centering
\includegraphics[width=1\linewidth]{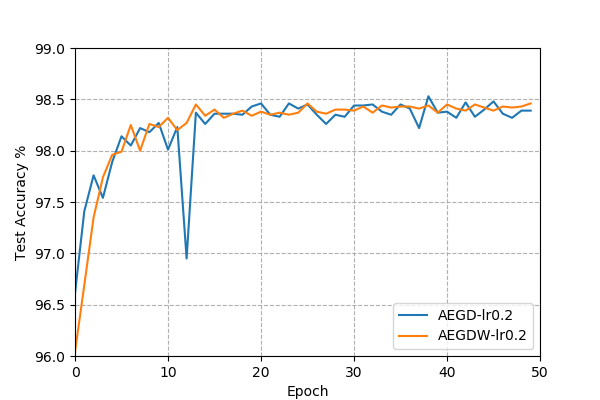}
\caption{AEGD v.s. AEGDW}
\end{subfigure}%
\begin{subfigure}[b]{0.5\linewidth}
\centering
\includegraphics[width=1\linewidth]{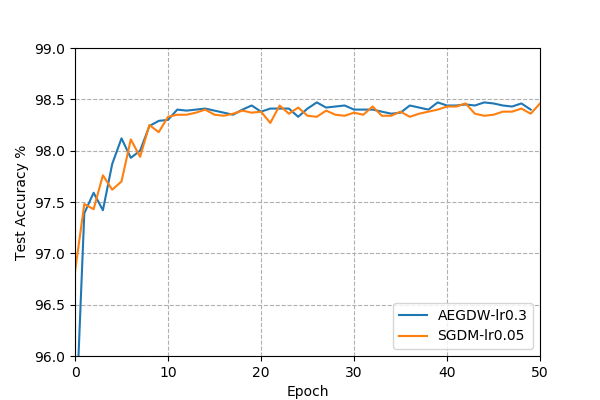}
\caption{Test Accuracy for MLP on MNIST}
\end{subfigure}%
\newline
\begin{subfigure}[b]{0.5\linewidth}
\centering
\includegraphics[width=1\linewidth]{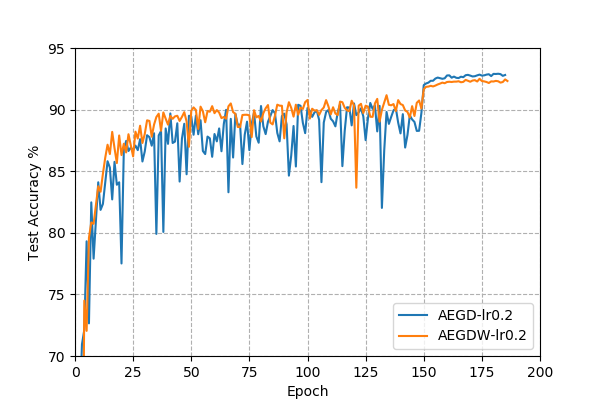}
\caption{AEGD v.s. AEGDW}
\end{subfigure}%
\begin{subfigure}[b]{0.5\linewidth}
\centering
\includegraphics[width=1\linewidth]{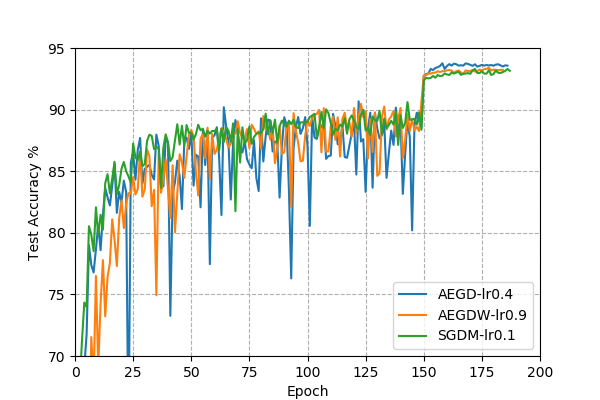}
\caption{Test Accuracy for CifarNet on CIFAR-10}
\end{subfigure}%
\captionsetup{format=hang}
\caption{Test accuracy for MLP on MNIST and CifarNet on CIFAR-10}
\label{fig:mnist_cifar}
\end{figure}

\textbf{MLP on MNIST.} We train a simple multi-layer perceptron (MLP), a special class of
feedforward neural networks,  with one hidden layer of $200$ neurons for the multi-class classification problem on MNIST data set. We run $50$ epochs with a batch size of $128$ and a weight decay of $10^{-4}$ for this  experiment. Figure \ref{fig:mnist_cifar} (b) shows that AEGDW performs slightly better than SGDM in this case. This is as expected for simple networks.

\textbf{CifarNet on CIFAR-10.}  We also train a simple 3-block convolutional neural network and name it CifarNet
on CIFAR-10. We use the same training set as before -- that is, reduce the learning rates by $10$ after $150$ epochs -- with a minibatch size of $128$ and weight decay of $10^{-4}$. Results for this experiment are reported in Figure \ref{fig:mnist_cifar} (d).  The overall performance of each algorithm for CifarNet on CIFAR-10 is similar to the experiments in Figure \ref{fig:cifar}.

AEGDW as an improved AEGD can help to reduce the variance and generalize better early, but such
generalization performance does not seem to sustain after decaying the learning rate as pre-scheduled (see examples in Figure \ref{fig:mnist_cifar} (a) and (c)).
 Overall, AEGDW or AEGD can give better generalization performance than SGDM as evidenced by our experiments.

\begin{algorithm}[t]
\caption{AEGD with decoupled weight decay (AEGDW). Good default setting for parameters are $c=1$ and $\eta=0.7/0.9$ for deep learning problems}
\label{alg:AEGDW}
\begin{algorithmic}[1] 
\Require $\{f_j({\theta})\}_{j=1}^m$, $\eta$: the step size, $\theta_0$: initial guess of $\theta$, and
$T$: the total number of iterations.
\Require $c$: a parameter such that $f({\theta})+c>0$ for all $\theta \in \mathbb{R}^n$, 
initial energy, $r_0=\sqrt{f(\theta_0)+c} {\bf 1}$,
weight decay factor $\lambda \in \mathbb{R}$
\For{$k=0$ to $k-1$}
\State $v_k:=\nabla f_{i_k}(\theta_k)/\big(2\sqrt{f_{i_k}(\theta_k)+c}\big)$ 
($i_k$ is a random sample from $[m]$ at step $k$)
\State $r_{k+1} = r_k/(1+2\eta v_k\odot v_k )$ (update energy)
\State $\theta_{k+1} = \theta_k - \eta (2r_{k+1}\odot v_k+\lambda\theta_k)$ (update parameters with weight decay)
\EndFor
\State \textbf{return} $\theta_k$
\end{algorithmic}
\end{algorithm}

\bibliographystyle{siamplain}
\bibliography{references}
\end{document}